\documentclass[a4paper,10pt]{article}

\def\expensiveFigures{1}
\def\longtitle{The~SCIP~Optimization~Suite~9.0}
\def\shortfunding{The work for this article has been partly conducted within the
  \emph{Research Campus MODAL} funded by the German Federal Ministry of
  Education and Research (BMBF grant number 05M14ZAM) and has received funding from
  the European Union's Horizon 2020 research and innovation programme under
  grant agreement No 773897. It has also been
  partly supported by
  the German Research Foundation (DFG) within the
  Collaborative Research Center 805, Project A4,
  and the EXPRESS project of the priority program CoSIP (DFG-SPP 1798),
  the German Research Foundation (DFG) within the project HPO-NAVI (project number 391087700).
  }

\makeatletter
\DeclareRobustCommand*{\escapeus}[1]{%
  \begingroup\@activeus\scantokens{#1 }\endgroup}
\begingroup\lccode`\~=`\_\relax
   \lowercase{\endgroup\def\@activeus{\catcode`\_=\active \let~\_}}
\makeatother

\usepackage[utf8]{inputenc}
\usepackage[font=small, labelfont=bf, margin=1cm]{caption}
\usepackage{amsmath,amsfonts,amssymb}
\usepackage{tabularx,booktabs}
\usepackage{url}
\usepackage{rotating}
\usepackage{microtype}

\usepackage[textsize=small,textwidth=3.3cm]{todonotes}

\usepackage{xspace}
\usepackage{etex}
\usepackage{enumitem}
\usepackage{dsfont}
\usepackage{mathtools}
\usepackage[linesnumbered,ruled,vlined]{algorithm2e}

\SetCommentSty{AlgoCommentStyle}
\usepackage{changes}
\usepackage{longtable}
\usepackage{tabularx}
\usepackage{varwidth}
\usepackage{listings}
\usepackage{diagbox}
\usepackage{makecell}

\usepackage{geometry}
\usepackage{amsthm}
\usepackage{multirow}

\usepackage[breaklinks,colorlinks=true,citecolor=black,linkcolor=black,urlcolor=black,hyperfootnotes=false,
  pdftitle={\longtitle},
  pdfsubject={\longtitle}]{hyperref}

\usepackage[capitalize]{cleveref}
\usepackage[numbers]{natbib}
\usepackage{doi}

\usepackage{tikz}
\usepackage{tikz-3dplot}
\usepackage{tikzscale}
\usetikzlibrary{arrows}
\usetikzlibrary{calc}
\usetikzlibrary{shapes}
\usetikzlibrary{snakes}
\usetikzlibrary{patterns}
\usetikzlibrary{backgrounds,topaths}
\usetikzlibrary{matrix,chains,scopes,positioning,arrows,fit}

\usepackage{pgfplots}
\pgfplotsset{compat=1.15}
\usepackage{subcaption}
\usepackage{textcomp}
\usepgfplotslibrary{groupplots}
\usepgfplotslibrary{units}

\usepackage{graphicx}

\makeatletter
\newcommand*\rel@kern[1]{\kern#1\dimexpr\macc@kerna}
\newcommand*\widebar[1]{%
  \begingroup
  \def\mathaccent##1##2{%
    \rel@kern{0.8}%
    \overline{\rel@kern{-0.8}\macc@nucleus\rel@kern{0.2}}%
    \rel@kern{-0.2}%
  }%
  \macc@depth\@ne
  \let\math@bgroup\@empty \let\math@egroup\macc@set@skewchar
  \mathsurround\z@ \frozen@everymath{\mathgroup\macc@group\relax}%
  \macc@set@skewchar\relax
  \let\mathaccentV\macc@nested@a
  \macc@nested@a\relax111{#1}%
  \endgroup
}
\makeatother

\newcommand{\LP}{{LP}\xspace}
\newcommand{\LPs}{{LPs}\xspace}

\newcommand{\CIPs}{{CIPs}\xspace}
\newcommand{\MIP}{{MIP}\xspace}
\newcommand{\MIPs}{{MIPs}\xspace}
\newcommand{\MILP}{{MILP}\xspace}
\newcommand{\MILPs}{{MILPs}\xspace}
\newcommand{\NLP}{{NLP}\xspace}

\newcommand{\MINLP}{{MINLP}\xspace}
\newcommand{\MINLPs}{{MINLPs}\xspace}

\newcommand{\MISDP}{{MISDP}\xspace}
\newcommand{\MISDPs}{{MISDPs}\xspace}

\newcommand{\T}{\top}

\newcommand{\abs}[1]{\lvert{#1}\rvert}

\newcommand{\defi}{\coloneqq}
\newcommand{\st}{\,:\,}

\newcommand{\iprod}[2]{\langle{#1},{#2}\rangle}

\newcommand{\linobj}{c}
\newcommand{\nonlinobj}{f}
\newcommand{\linmatrix}{A}
\newcommand{\nonlincons}{g}
\newcommand{\rhs}{b}
\newcommand{\lb}{\ell}
\newcommand{\ub}{u}

\newcommand{\consindex}{\mathcal{M}}
\newcommand{\varindex}{\mathcal{N}}
\newcommand{\intvarindex}{\mathcal{I}}

\newcommand{\R}{\mathds{R}}

\newcommand{\Z}{\mathds{Z}}
\newcommand{\Rinf}{\ensuremath{\widebar{\mathds{R}}}\xspace}

\usepackage{siunitx}
\newcommand{\cleaninst}{all}
\newcommand{\affected}{affected}
\newcommand{\alloptimal}{{both-solved}\xspace}
\newcommand{\difftimeouts}{{diff-timeouts}\xspace}
\newcommand{\nonconvex}{nonconvex}
\newcommand{\convex}{convex}
\newcommand{\bracket}[2]{[#1,#2]}

\usepackage{xparse}
\ExplSyntaxOn
\def\myround#1{\num{\fp_eval:n {round(#1, 2)}}}
\ExplSyntaxOff

\definecolor{c1}{HTML}{000060}
\definecolor{c2}{HTML}{0000FF}
\definecolor{c3}{HTML}{36648B}
\definecolor{c4}{HTML}{4682B4}
\definecolor{c5}{HTML}{5CACEE}
\definecolor{c6}{HTML}{00FFFF}
\definecolor{c7}{HTML}{008888}
\definecolor{c8}{HTML}{00DD99}
\definecolor{c9}{HTML}{527B10}
\definecolor{c10}{HTML}{7BC618}
\definecolor{c11}{HTML}{33AA00}

\definecolor{scipoldcol}{HTML}{36648B}
\definecolor{scipnewcol}{HTML}{7BC618}

\newcommand{\solver}[1]{\textsc{#1}\xspace}
\newcommand{\scipopt}{\scip Optimization Suite\xspace}
\newcommand{\scipversion}{9.0}

\newcommand{\scipoptv}{\scipopt~\scipversion\xspace}
\newcommand{\scip}{\solver{SCIP}}
\newcommand{\scipv}{\solver{SCIP}~\scipversion\xspace}
\newcommand{\soplex}{\solver{SoPlex}}
\newcommand{\soplexversion}{6.0}
\newcommand{\soplexv}{\solver{SoPlex}~\soplexversion\xspace}
\newcommand{\papilo}{\solver{PaPILO}}
\newcommand{\papiloversion}{2.2}
\newcommand{\papilov}{\solver{PaPILO}~\papiloversion\xspace}
\newcommand{\zimpl}{\solver{Zimpl}} 

\newcommand{\ug}{\solver{UG}}
\newcommand{\presollib}{\solver{PaPILO}}

\newcommand{\gcg}{\solver{GCG}}
\newcommand{\gcgversion}{3.6}
\newcommand{\gcgv}{\gcg~\gcgversion\xspace}

\newcommand{\scipsdp}{\solver{SCIP-SDP}}

\newcommand{\scipjack}{\solver{SCIP-Jack}}

\newcommand{\scipjl}{\solver{SCIP.jl}}

\newcommand{\pyscipopt}{\solver{PySCIPOpt}}
\newcommand{\pyscipoptml}{\solver{PySCIPOpt-ML}}
\newcommand{\russcip}{\solver{russcip}} 
\newcommand{\scippp}{\solver{SCIP++}} 
\newcommand{\exscip}{\solver{ExactSCIP}} 
\newcommand{\jscipopt}{\solver{JSCIPOpt}}
\newcommand{\pysoplex}{\solver{PySoPlex}}
\newcommand{\fscip}{\solver{FiberSCIP}}
\newcommand{\pscip}{\solver{ParaSCIP}}

\newcommand{\vipr}{\solver{Vipr}}

\newcommand{\ipopt}{\solver{Ipopt}}
\newcommand{\cppad}{\solver{CppAD}}

\newcommand{\highs}{\solver{HiGHS}}

\newcommand{\veripb}{\solver{VeriPB}}

\newcommand{\nbsc}[1]{\mbox{#1}\xspace}
\newcommand{\miplib}{\nbsc{MIPLIB}}

\newcommand{\minlplibtwo}{\nbsc{MINLPLib}}

\definecolor{darkgreen}{HTML}{008800}

\newcommand{\fa}{\text{ for all }}

\newcommand{\bliss}{\solver{bliss}}

\newcommand{\nauty}{\solver{nauty}}
\newcommand{\sassy}{\solver{sassy}}

\newcommand{\sign}{\mathrm{sgn}}

\theoremstyle{plain}

\newtheorem{definition}{Definition}[subsection]

\setlist[itemize]{leftmargin=3.45ex}
\setlist[itemize,1]{label=$-$,itemsep=0ex,topsep=0.9ex}
\setlist[itemize,2]{label=$\cdot$,topsep=0.5ex,leftmargin=2.75ex}
\setlist[enumerate]{leftmargin=3ex,itemsep=0.1ex,parsep=1ex,topsep=0.9ex}

\definecolor{tabcolor}{HTML}{6666AA}
\definecolor{f1}{HTML}{000060}
\definecolor{f2}{HTML}{0000FF}
\definecolor{f3}{HTML}{36648B}
\definecolor{f4}{HTML}{4682B4}
\definecolor{f5}{HTML}{5CACEE}
\definecolor{f6}{HTML}{00FFFF}
\definecolor{f7}{HTML}{00DD99}
\definecolor{f8}{HTML}{008888}
\definecolor{f9}{HTML}{000000}

\usepackage{float}
\newfloat{program}{thp}{lop}
\floatname{program}{Program}
\crefname{program}{program}{programs}

\newcommand{\inputExpensiveFigure}[1]{
\ifthenelse{\expensiveFigures = 1}{\input{#1}}{}
}

\usepackage{titling}
\pretitle{\begin{center}\linespread{1.05}\Large}
\posttitle{\par\end{center}\vskip 0.5em}
\preauthor{\begin{center}\linespread{1.15}\large \lineskip 0.6em\begin{tabular}[t]{c}}
\postauthor{\end{tabular}\par\end{center}\vskip 0.6em}
\predate{\begin{center}}
\postdate{\par\end{center}}
\usepackage[]{titlesec}
\usepackage{etoolbox}
\makeatletter
\patchcmd{\ttlh@hang}{\parindent\z@}{\parindent\z@\leavevmode}{}{}
\patchcmd{\ttlh@hang}{\noindent}{}{}{}
\makeatother
\titleformat{\section}
{\normalfont\large\bfseries}{\thesection}{0.9ex}{}
\titleformat{\subsection}
{\normalfont\normalsize\bfseries}{\thesubsection}{0.9ex}{}
\titleformat{\subsubsection}
{\normalfont\normalsize\upshape}{\thesubsubsection}{0.9ex}{}
\titleformat{\paragraph}[runin]
{\normalfont\normalsize\itshape}{\theparagraph}{1em}{}
\titleformat{\subparagraph}[runin]
{\normalfont\normalsize\itshape}{\theparagraph}{1em}{}
\titlespacing*{\section}     {0pt}{21dd plus 8pt minus 4pt}{10.5dd}
\titlespacing*{\subsection}   {0pt}{21dd plus 8pt minus 4pt}{10.5dd}
\titlespacing*{\subsubsection}{0pt}{19dd plus 8pt minus 4pt}{10.5dd}
\titlespacing*{\paragraph}   {0pt}{13pt plus 8pt minus 4pt}{1em}
\titlespacing*{\subparagraph}   {0pt}{13pt plus 8pt minus 4pt}{1em}

\usepackage[page]{appendix} 

\newcommand{\myand}{$\cdot$\xspace}
\newcommand{\myorcidlink}[1]{\,\href{https://orcid.org/#1}{\raisebox{-0.45ex}{\includegraphics[width=1.8ex]{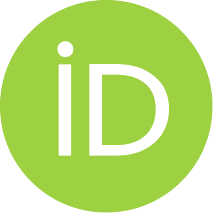}}}}

\begin{document}

\title{\longtitle}

\author{%
  Suresh Bolusani\protect\myorcidlink{0000-0002-5735-3443} \myand
  Mathieu Besançon\protect\myorcidlink{0000-0002-6284-3033} \myand
  Ksenia Bestuzheva\protect\myorcidlink{0000-0002-7018-7099} \and
  Antonia Chmiela\protect\myorcidlink{0000-0002-4809-2958} \myand
  Jo{\~{a}}o~Dion{\'{i}}sio\protect\myorcidlink{0009-0005-5160-0203} \myand
  Tim Donkiewicz\protect\myorcidlink{0000-0002-5721-3563} \and
  Jasper van Doornmalen\protect\myorcidlink{0000-0002-2494-0705} \myand
  Leon Eifler\protect\myorcidlink{0000-0003-0245-9344} \myand
  Mohammed Ghannam\protect\myorcidlink{0000-0001-9422-7916} \and
  Ambros Gleixner\protect\myorcidlink{0000-0003-0391-5903} \myand
  Christoph Graczyk\protect\myorcidlink{0000-0001-8990-9912} \myand
  Katrin Halbig\protect\myorcidlink{0000-0002-8730-3447} \and
  Ivo Hedtke\protect\myorcidlink{0000-0003-0335-7825} \myand
  Alexander Hoen\protect\myorcidlink{0000-0003-1065-1651} \myand
  Christopher Hojny\protect\myorcidlink{0000-0002-5324-8996} \and
  Rolf van der Hulst\protect\myorcidlink{0000-0002-5941-3016} \myand
  Dominik Kamp\protect\myorcidlink{0009-0005-5577-9992} \myand
  Thorsten Koch\protect\myorcidlink{0000-0002-1967-0077} \and
  Kevin Kofler \myand
  Jurgen Lentz\protect\myorcidlink{0009-0000-0531-412X} \myand
  Julian Manns \myand
  Gioni Mexi\protect\myorcidlink{0000-0003-0964-9802} \and
  Erik M\"uhmer\protect\myorcidlink{0000-0003-1114-3800} \myand
  Marc E. Pfetsch\protect\myorcidlink{0000-0002-0947-7193} \myand
  Franziska Schlösser \and
  Felipe Serrano\protect\myorcidlink{0000-0002-7892-3951} \myand
  Yuji Shinano\protect\myorcidlink{0000-0002-2902-882X} \myand
  Mark Turner\protect\myorcidlink{0000-0001-7270-1496} \and
  Stefan Vigerske\protect\myorcidlink{0009-0001-2262-0601} \myand
  Dieter Weninger\protect\myorcidlink{0000-0002-1333-8591}  \myand
  Liding Xu\protect\myorcidlink{0000-0002-0286-1109}
  \thanks{Extended author information is available at the end of the paper.
    \shortfunding}}

\date{26 February 2024}

\newgeometry{left=38mm,right=38mm,top=35mm}

\maketitle

\paragraph{\bf Abstract}

The \scipopt provides a collection of software packages for mathematical optimization, centered around the constraint integer programming (CIP) framework \scip.
This report discusses the enhancements and extensions included in \scipoptv.
The updates in \scipv include improved symmetry handling, additions and improvements of nonlinear handlers and primal heuristics, a new cut generator and two new cut selection schemes, a new branching rule, a new \LP interface, and several bugfixes.
\scipoptv also features new Rust and C++ interfaces for \scip, new Python interface for \soplex, along with enhancements to existing interfaces.
\scipoptv also includes new and improved features in the \LP solver \soplex, the presolving library \papilo, the parallel framework \ug, the decomposition framework \gcg, and the \scip extension \scipsdp.
These additions and enhancements have resulted in an overall performance improvement of \scip in terms of solving time, number of nodes in the branch-and-bound tree, as well as the reliability of the solver.



\newpage

\section{Introduction}
\label{sect:introduction}

The \scipopt comprises a set of complementary software packages designed to model and
solve a large variety of mathematical optimization problems:
\begin{itemize}
\item the constraint integer programming solver
\scip~\cite{Achterberg2009}, a solver for
mixed-integer linear and nonlinear programs as well as a flexible framework for
branch-cut-and-price,
\item the simplex-based linear programming solver
  \soplex~\cite{Wunderling1996},
\item the modeling language \zimpl~\cite{Koch2004},
\item the presolving library \presollib for linear and mixed-integer linear programs,
\item the automatic decomposition solver \gcg~\cite{GamrathLuebbecke2010}, and
\item the \ug framework for parallelization of branch-and-bound
  solvers~\cite{Shinano2018}.
\end{itemize}

All six tools are freely available as open-source software packages, either using the Apache~2.0 or the GNU Lesser General Public License.
There also exist two notable continuously developed extensions to the \scipopt: the award-winning Steiner tree solver
\scipjack~\cite{Gamrath2017scipjack} and the mixed-integer semidefinite
programming solver \scipsdp~\cite{GallyPfetschUlbrich2018}.
This report describes the improvements and new features contained in version~9.0 of the \scipopt.

\paragraph{Background}

\scip is designed as a solver for \emph{constraint integer programs} (\CIPs), a generalization of mixed-integer linear and nonlinear programs (\MILPs and \MINLPs). \CIPs are finite-dimensional optimization problems with arbitrary constraints and a linear objective function that satisfy the following property: if all integer variables are fixed, the remaining subproblem must form a linear or nonlinear program (\LP or \NLP).
To solve \CIPs, \scip constructs relaxations---typically linear relaxations, but also nonlinear relaxations are possible, or relaxations based on semidefinite programming for \scipsdp. If the relaxation solution is not feasible for the current subproblem, an \emph{enforcement} procedure is called that takes measures to resolve the infeasibility, for example by branching or by separating cutting planes.

The most important subclass of \CIPs that are solvable with \scip are \emph{mixed-integer programs} (\MIPs) which can be purely linear (\MILPs) or contain nonlinearities (\MINLPs).
\MILPs are optimization problems of the form
\begin{equation}
  \begin{aligned}
    \min \quad& \linobj^\T x \\
    \text{s.t.} \quad& \linmatrix x \geq \rhs, \\
    &\lb_{i} \leq x_{i} \leq \ub_{i} && \fa i \in \varindex, \\
    &x_{i} \in \Z && \fa i \in \intvarindex,
  \end{aligned}
  \label{eq:generalmilp}
\end{equation}
defined by $c \in \R^n$, $A \in\R^{m\times n}$, $ \rhs\in \R^{m}$, $\lb$, $\ub \in
\Rinf^{n}$, and the index set of integer variables $\mathcal{I} \subseteq \mathcal{N} \defi \{1, \ldots, n\}$.  The usage of $\Rinf \defi \R \cup
\{-\infty,\infty\}$ allows for variables that are free or bounded only in
one direction (we assume that variables are not fixed to~$\pm \infty$).
In contrast, \MINLPs are optimization problems of the form
\begin{equation}
  \begin{aligned}
    \min \quad& \nonlinobj(x) \\
    \text{s.t.} \quad& \nonlincons_{k}(x) \leq 0 && \fa k \in \consindex, \\
    &\lb_{i} \leq x_{i} \leq \ub_{i} && \fa i \in \varindex, \\
    &x_{i} \in \Z && \fa i \in \intvarindex,
  \end{aligned}
  \label{eq:generalminlp}
\end{equation}
where the functions $\nonlinobj : \R^n \rightarrow \R$ and $\nonlincons_{k} : \R^{n}
\rightarrow \R$, $k \in \consindex \defi \{1,\ldots,m\}$, are possibly nonconvex.
Within \scip, we assume that $\nonlinobj$ is linear and that $\nonlincons_k$ are specified explicitly in
algebraic form using a known set of base expressions.

Due to its design as a solver for \CIPs, \scip can be extended by plugins for more general or problem-specific classes of optimization problems.
The core of \scip is formed by a central branch-cut-and-price algorithm that utilizes an \LP as the default relaxation which can be solved by a number of different \LP solvers, controlled through a uniform \emph{\LP interface}. To be able to handle any type of constraint, a \emph{constraint handler} interface is provided, which allows for the integration of new constraint types, and provides support for many different well-known types of constraints out of the box.
Further, advanced solving methods like primal heuristics, branching rules, and cutting plane separators can also be integrated as plugins with a pre-defined interface.
\scip comes with many such plugins needed to achieve a good \MILP and \MINLP performance.
In addition to plugins supplied as part of the SCIP distribution, new plugins can be created by users.
The design approach and solving process is described in detail by
Achterberg~\cite{Achterberg2007a}.

Although it is a standalone solver,
\scip interacts closely with the other components of the \scipopt.
\zimpl is integrated into \scip as a reader plugin, making it possible to read \zimpl problem instances directly by \scip.
\presollib is integrated into \scip as an additional presolver plugin.
The \LPs that need to be solved as relaxations in the branch-and-bound process are by default solved with \soplex.
Interfaces to most actively developed external \LP solvers exist, and new interfaces can be added by users.
\gcg extends \scip to automatically detect problem structure
and generically apply decomposition algorithms based on the Dantzig-Wolfe or the
Benders' decomposition schemes.
Finally, the default instantiations of the \ug framework use \scip as a base
solver in order to perform branch-and-bound in parallel computing
environments with shared or distributed memory architectures.

\paragraph{New Developments and Structure of the Report}
This report is structured into three main parts.
First, the changes and progress made in the solving process of \scip are explained
and the resulting performance improvements on \MILP and \MINLP instances are analyzed, both in terms of performance and robustness.
A performance comparison of \scipv against \scip 8.0 is carried out in Section \ref{sect:performance}.

Second, improvements to the core of \scip are presented in Section~\ref{sect:scip}, which include
\begin{itemize}
   \item improved symmetry handling on non-binary variables, symmetry handling for custom constraints, and signed permutation symmetries,
   \item symmetry preprocessing using the new interfaces to \nauty and \sassy,
   \item a new constraint handler for signomial inequalities as well as cut-strengthening for quadratic expressions,
   \item a new indicator diving heuristic, extensions to the existing dynamic partition search heuristic, as well as a new online scheduling feature for primal heuristics,
   \item a new Lagromory separator, as well as improvements in cut selection,
   \item a new branching criterion called GMI branching that is incorporated into the existing scoring function and acts as a tie-breaker for the existing branching rules,
   \item a new interface to the \highs~ \LP solver \cite{HuangfuHall15}, and
   \item certain technical improvements.
\end{itemize}

Third, improvements to the other components of the \scipopt and extensions to the interfaces are presented.
Improvements to the default \LP solver \soplex and presolver \papilo are explained in Sections \ref{sect:soplex}
and \ref{sect:papilo}, respectively. Extensions to the interfaces of \scip are presented in Section \ref{sect:interfaces}. Besides improvements and extensions to existing interfaces, this section includes two new interfaces for \scip: (1) \russcip~\cite{russcipgithub}, a new Rust interface, and (2) \scippp~\cite{scipppgithub}, a new C++ interface, as well as a new Python interface for \soplex called \pysoplex~\cite{pysoplexgithub}.
Improvements to distributed computing with \ug and to Dantzig-Wolfe decompositions with \gcg are presented in Sections \ref{sect:ug}, and \ref{sect:gcg}, respectively;
and finally updates to the \scip extension \scipsdp for semidefinite problems are presented in Section \ref{sect:scip-sdp}.
Not included in this release, but available as a beta version, is \exscip~\cite{exscipgithub}, a new extension of \scip that allows for the exact solution of MILPs with rational input data without roundoff errors and zero numerical tolerances.

\section{Overall Performance Improvements for MILP and MINLP}
\label{sect:performance}
In this section, we present computational experiments conducted by running \scip without parameter tuning or
algorithmic variations to assess the performance changes since the 8.0.0 release.
We detail below the methodology and results of these experiments.

The indicators of interest to compare the two versions of \scip
on a given subset of instances are
the number of solved instances, the shifted geometric mean
of the number of branch-and-bound nodes, and the shifted geometric mean
of the solving time.
The \emph{shifted geometric mean} of values $t_1, \dots, t_n$ is
\[
\left(\prod_{i=1}^n(t_i + s)\right)^{1/n} - s.
\]
The shift~$s$ is set to 100~nodes and 1~second, respectively.

\subsection{Experimental Setup}

As baseline we use \scip~8.0.0, with \soplex 6.0.0 as the underlying
\LP solver, and \papilo 2.0.0 for enhanced presolving.
We compare it with \scip~9.0.0 with \soplex 7.0.0 and \papilo 2.2.0.
Both \scip versions were compiled using \gcg 10.2.1, use \ipopt 3.14.14 as NLP subsolver built with \solver{HSL MA27} as linear system solver, \solver{Intel MKL} as linear algebra package,
\cppad 20180000.0 as algorithmic differentiation library,
and \bliss 0.77 for graph automorphisms to detect symmetry in \MIPs.
\scip~9.0.0 additionally uses \sassy 1.1 as a preprocessor for \bliss.
The time limit was set to 7200 seconds in all cases.
Furthermore, for \MINLP, a relative gap limit of $10^{-4}$ and an absolute gap limit of $10^{-6}$ were set.

The \MILP instances were selected from \miplib 2017~\cite{miplib2017}, including all instances
previously solved by previous \scip versions with at least one of five random seeds or newly solved by \scip~9.0.0 with at least one of five random seeds;
this amounted to 158~instances.
The \MINLP instances were selected in a similar way from the \minlplibtwo\footnote{\url{https://www.minlplib.org}} for a total of 179 instances.

All performance runs were carried out on identical machines with Intel Xeon Gold 5122 @ 3.60GHz
and 96GB RAM. A single run was carried out on each machine in a single-threaded mode.
Each instance was solved with \scip using five different seeds for random number generators.
This results in a testset of 790 \MILPs and 810 \MINLPs.
Instances for which the solver reported numerically inconsistent results are excluded from the results below.

\subsection{\MILP Performance}
Table~\ref{tbl:MIP5seeds} presents a comparison between \scipv and \scip~8.0 regarding their \MILP performance.
\scipv improves the solving capabilities for \MILP by solving 19 more instances than \scip~8.0.  In terms of the shifted geometric mean of the  running time, both versions perform almost equally across all instances, with \scipv being 2\% faster on affected instances.
On the subset of harder instances in the \bracket{1000}{7200} bracket, i.e., instances that take at least $1000$ seconds to be solved with at least one setting, the speedup is larger and amounts to 6\%. 
To compare average tree size across the two versions, we restrict to the \alloptimal subset since the number of nodes for instances that time out is not easy to interpret. In the \alloptimal subset, \scipv significantly reduces the average tree size by 17\%.
Finally, it's worth noting that \scipv incorporates a large number of bugfixes. While these bugfixes have introduced a slowdown, they contribute significantly to the overall reliability of the solver.

\begin{table}[h]
    \caption{Performance comparison of \scip~9.0 and \scip~8.0 for \MILP instances}
    \label{tbl:MIP5seeds}
    \scriptsize
    
    \begin{tabular*}{\textwidth}{@{}l@{\;\;\extracolsep{\fill}}rrrrrrrrr@{}}
    \toprule
    &           & \multicolumn{3}{c}{\scip~9.0.0+\soplex~7.0.0} & \multicolumn{3}{c}{\scip~8.0.0+\soplex~6.0.0} & \multicolumn{2}{c}{relative} \\
    \cmidrule{3-5} \cmidrule{6-8} \cmidrule{9-10}
    Subset                & instances &                                   solved &       time &        nodes &                                   solved &       time &        nodes &       time &        nodes \\
    \midrule
    \cleaninst            &       785 &                                      637 &      433.9 &         4307 &                                      618 &      439.0 &         5236 &        1.01 &          1.22 \\
    \affected             &       647 &                                      610 &      297.1 &         3874 &                                      591 &      301.8 &         4763 &        1.02 &          1.23 \\
    \cmidrule{1-10}
    \bracket{0}{tilim}    &       674 &                                      637 &      272.9 &         3332 &                                      618 &      276.7 &         4065 &        1.01 &          1.22 \\
    \bracket{1}{tilim}    &       669 &                                      632 &      283.8 &         3423 &                                      613 &      287.7 &         4182 &        1.01 &          1.22 \\
    \bracket{10}{tilim}   &       617 &                                      580 &      399.8 &         4463 &                                      561 &      404.9 &         5567 &        1.01 &          1.25 \\
    \bracket{100}{tilim}  &       460 &                                      423 &      986.9 &        11282 &                                      404 &      985.4 &        14244 &        1.00 &          1.26 \\
    \bracket{1000}{tilim} &       278 &                                      241 &     2240.0 &        30971 &                                      222 &     2383.3 &        41101 &        1.06 &          1.33 \\
    \difftimeouts         &        93 &                                       56 &     3899.0 &       100922 &                                       37 &     4592.9 &       160892 &        1.18 &          1.59 \\
    \alloptimal           &       581 &                                      581 &      178.1 &         1897 &                                      581 &      176.2 &         2220 &        0.99 &          1.17 \\
            \bottomrule
    \end{tabular*}
    \end{table}
\subsection{\MINLP Performance}
Table~\ref{tbl:MINLP5seeds} summarizes the results for the performance of \scipv as compared to \scip8.0 for the \MINLP instances. Besides increasing the number of solved instances by $5$, the changes introduced in \scipv improve the performance of \scip in both the overall solving time as well as the number of nodes needed. On the whole testset, \scipv improves the performance by about $4 \%$ in time and by $13 \%$ in nodes, both in shifted geometric means. This improvement increases with the difficulty of the instances: When looking at the most difficult testset $\bracket{1000}{7200}$, \scipv outperforms \scip8.0 by $20 \%$ and $46 \%$ in the solving time and nodes in shifted geometric means, respectively. Furthermore, the improvement in solving time is mainly observed on nonconvex instances. In particular, \scipv is $8 \%$ faster than \scip8.0 on nonconvex instances whereas both versions perform almost equally when only convex instances are considered.

\begin{table}[t]
\caption{Performance comparison of \scip~9.0 and \scip~8.0 for \MINLP instances}
\label{tbl:MINLP5seeds}
\scriptsize

\begin{tabular*}{\textwidth}{@{}l@{\;\;\extracolsep{\fill}}rrrrrrrrr@{}}
\toprule
&           & \multicolumn{3}{c}{\scip~8.0.0+\soplex~6.0.0} & \multicolumn{3}{c}{\scip~9.0.0+\soplex~7.0.0} & \multicolumn{2}{c}{relative} \\
\cmidrule{3-5} \cmidrule{6-8} \cmidrule{9-10}
Subset                & instances &                                   solved &       time &        nodes &                                   solved &       time &        nodes &       time &        nodes \\
\midrule
\cleaninst            &       839 &                                      810 &       32.6 &         2800 &                                      815 &       31.2 &         2489 &        1.04 &          1.13 \\
\affected             &       783 &                                      770 &       32.4 &         2949 &                                      775 &       31.4 &         2608 &        1.03 &          1.13 \\
\cmidrule{1-10}
\bracket{0}{7200}    &       823 &                                      810 &       29.3 &         2613 &                                      815 &       28.0 &         2323 &        1.04 &          1.12 \\
\bracket{1}{7200}    &       767 &                                      754 &       36.6 &         3288 &                                      759 &       34.9 &         2907 &        1.05 &          1.13 \\
\bracket{10}{7200}   &       529 &                                      516 &       97.0 &         6333 &                                      521 &       90.2 &         5498 &        1.08 &          1.15 \\
\bracket{100}{7200}  &       248 &                                      235 &      489.1 &        29290 &                                      240 &      414.0 &        21767 &        1.18 &          1.35 \\
\bracket{1000}{7200} &        96 &                                       83 &     2020.8 &       110115 &                                       88 &     1682.8 &        75556 &        1.20 &          1.46 \\
\difftimeouts         &        21 &                                        8 &     3492.6 &        72059 &                                       13 &     1691.3 &        21896 &        2.07 &          3.29 \\
\alloptimal           &       802 &                                      802 &       25.7 &         2389 &                                      802 &       25.1 &         2187 &        1.03 &          1.09 \\
\cmidrule{1-10}
\convex               &       168 &                                      163 &       31.4 &         3601 &                                      165 &       31.0 &         3177 &        1.01 &          1.13 \\
\nonconvex            &       571 &                                      547 &       38.0 &         2783 &                                      550 &       35.2 &         2443 &        1.08 &          1.14 \\
\bottomrule
\end{tabular*}
\end{table}

\section{SCIP}
\label{sect:scip}

\subsection{Symmetry Handling}
\label{sec:symmetry}
Symmetries of an \MILP or \MINLP are maps that transform
feasible solutions into feasible solutions with the same objective value.
When not handled appropriately, such symmetries deteriorate the
performance of (spatial) branch-and-bound algorithms since symmetric
solutions are found and symmetric subproblems are explored repeatedly
without providing additional information to the solver.
The previous versions of \scip have already contained many state-of-the-art
algorithms to handle symmetries of binary variables and some basic cutting
planes to also handle symmetries of integer or continuous variables.

\scipv substantially extends the ability to handle symmetries
in three directions.
First, more sophisticated techniques are available to handle symmetries of
non-binary variables.
Second, the mechanism to detect symmetries has been completely
restructured.
While previous versions of \scip could only detect symmetries of the
available classes of constraints, \scipv can
also detect symmetries of custom constraints added by users.
Moreover, \scipv can also detect so-called signed permutation
symmetries, whereas previous versions could only detect permutation
symmetries.
Finally, to detect symmetries, \scip makes use of external software for
detecting graph automorphisms.
In the latest version, new interfaces to \nauty~\cite{Nauty}
as well as the preprocessor \sassy~\cite{Sassy} have been added.

\subsubsection{Symmetry Handling Methods}
\label{sec:symmetryhandling}
For the ease of exposition, consider an \MILP~$\min\,\{\linobj^\T x
\st Ax \geq b,\; x \in \Z^p \times \R^{n-p}\}$, where~$A \in \R^{m \times
  n}$, $b \in \R^m$, and~$c \in \R^n$.
Symmetries in \MILP are commonly permuting variables, i.e.,
a permutation~$\gamma$ of~$[n]$ acts on~$x \in \R^n$
as~$\gamma(x) = (x_{\gamma^{-1}(1)},\dots,x_{\gamma^{-1}(n)})$.
Permutation~$\gamma$ is called a \emph{formulation symmetry} of the \MILP if there
is a permutation~$\pi$ of~$[m]$ such that~$\gamma(c) = c$, $\pi(b) =
b$, $A_{\pi^{-1}(i),\gamma^{-1}(j)} = A_{i,j}$ for all~$(i,j) \in [m]
\times [n]$, and~$\gamma(i) \in [p]$ for all~$i \in [p]$.
Formulation symmetries can be detected by constructing a colored graph
whose color-preserving automorphisms correspond to symmetries of the
\MILP, see~\cite{Salvagnin2005} and Section~\ref{sec:symdetect} for more details.
Moreover, the definition of formulation symmetries can be extended to
\MINLPs by keeping the representation of
nonlinear constraints via expression trees invariant~\cite{Liberti2012}.

In previous versions of \scip, three main classes of symmetry handling
methods have been available:
\begin{enumerate}
\item Propagation and separation algorithms for
  the symmetry handling constraints orbisack~\cite{Hojny2020,HojnyPfetsch2019,KaibelLoos2011,Loos2010},
  symresack~\cite{Hojny2020,HojnyPfetsch2019},
  and orbitope~\cite{BendottiEtAl2021,KaibelPfetsch2008,KaibelEtAl2011};
  these constraints have only been able to handle symmetries of binary variables and
  enforce symmetry reductions based on a scheme that is
  determined before the solving process starts.

\item The propagation method orbital
  fixing~\cite{Margot2003,OstrowskiEtAl2011} to handle symmetries of
  binary variables; the corresponding symmetry reduction scheme is
  determined dynamically during the solving process.

\item Schreier-Sims table cuts (SST
  cuts)~\cite{LibertiOstrowski2014,Salvagnin2018}, which are cutting planes
  that are added to the \MILP/\MINLP and can handle symmetries of arbitrary
  variable types.
\end{enumerate}
Note that the first two classes are not compatible with each other due
to the different symmetry handling schemes.

\scipv features an implementation of a generalization of the first two
classes of methods as discussed in~\cite{vDoornmalenHojny2022}.
This generalization allows to also handle symmetries of non-binary
variables and to apply both classes simultaneously.
At the time of adding this new feature, the performance of \scip improved
by \SI{5.9}{\percent} on the \miplib 2017 benchmark test set; on the hard
instances that take at least \SI{1000}{\second} to be solved, the running
time even improved by \SI{25.4}{\percent}.

\subsubsection{Symmetry Detection}
\label{sec:symdetect}

In previous versions of \scip, symmetries could only be detected if all
types of constraints present in an \MILP or \MINLP are known by \scip, i.e.,
constraints whose corresponding constraint handler is part of the \scip
release.
In particular, this means that no symmetries could be detected and
automatically handled by \scip in the presence of custom constraints.
For \scipv, the symmetry detection mechanism has been
restructured.
Constraint handlers now support two optional callbacks
\texttt{CONSGETPERMSYMGRAPH} and \texttt{CONSGETSIGNEDPERMSYMGRAPH} that
allow constraint handlers to inform \scip about symmetries of their
constraints.
The former is used for the detection of permutation symmetries, whereas the
latter allows to detect signed permutation symmetries that we define
below.
During run time, \scip checks whether the constraint handlers of all
constraints present in a problem implement the new callbacks.
If this is the case, symmetries are detected; otherwise, symmetry detection
is disabled.

\paragraph{Detection of Permutation Symmetries}
As briefly explained above, permutation symmetries of an \MILP or \MINLP can
be detected by finding automorphisms of a suitable colored graph, which we
call the symmetry detection graph.
In the following, we explain this mechanism and how it can be implemented
using the \texttt{CONSGETPERMSYMGRAPH} callback.
We illustrate the ideas using the simple \MILP
$\max\,\{y + z \st -2w + 2x + 3y + 3z \leq 4,\; y,\, z \in \{0,1\} \}$.

To detect symmetries, every constraint defines its own local symmetry
detection graph.
Such a graph contains a colored node for every variable that is present in
the constraint as well as further colored nodes and edges that model
dependencies between the different variables.
The graph for a constraint qualifies as a symmetry detection graph for
\scip if it is connected and the restriction of every color-preserving
automorphism to the variable nodes corresponds to a permutation symmetry of
the corresponding constraint.
Moreover, two symmetry detection graphs are only allowed to be isomorphic
if their constraints are equivalent.

\begin{figure}[t]
  \centering
  \begin{subfigure}{0.3\textwidth}
    \centering
    \begin{tikzpicture}
      [n/.style={inner sep=0mm,minimum size=2mm,draw=black,circle}]
      \node (w) at (0,0) [n,label=left:{$w$},fill=red] {};
      \node (x) at (1,0) [n,label=left:{$x$},fill=blue] {};
      \node (y) at (2,0) [n,label=left:{$y$},fill=orange] {};
      \node (z) at (3,0) [n,label=left:{$z$},fill=orange] {};
      \node (rhs) at (1.5,1) [n,label=left:{rhs},fill=black] {};
      \node (c) at (1.5,2) [n,label=left:{cons}] {};

      \draw[-,draw=cyan,thick] (w) -- (rhs);
      \draw[-,draw=purple,thick] (x) -- (rhs);
      \draw[-,draw=green,thick] (y) -- (rhs);
      \draw[-,draw=green,thick] (z) -- (rhs);
      \draw[-,thick] (c) -- (rhs);
    \end{tikzpicture}
    \caption{permutation symmetries}
    \label{fig:detectpermsym}
  \end{subfigure}
  \hfill
  \begin{subfigure}{0.65\textwidth}
    \centering
    \begin{tikzpicture}
      [n/.style={inner sep=0mm,minimum size=2mm,draw=black,circle}]
      \node (w) at (0,0) [n,label=below:{$w$},fill=red] {};
      \node (w2) at (1,0) [n,label=below:{$-w$},fill=blue] {};
      \node (x) at (2,0) [n,label=below:{$x$},fill=blue] {};
      \node (x2) at (3,0) [n,label=below:{$-x$},fill=red] {};
      \node (y) at (4,0) [n,label=below:{$y$},fill=orange] {};
      \node (y2) at (5,0) [n,label=below:{$-y$},fill=black!50] {};
      \node (z) at (6,0) [n,label=below:{$z$},fill=orange] {};
      \node (z2) at (7,0) [n,label=below:{$-z$},fill=black!50] {};
      \node (rhs) at (3.5,1) [n,label=left:{rhs},fill=black] {};
      \node (c) at (3.5,1.8) [n,label=left:{cons}] {};

      \draw[-,draw=cyan,thick] (w) -- (rhs);
      \draw[-,draw=purple,thick] (x) -- (rhs);
      \draw[-,draw=green,thick] (y) -- (rhs);
      \draw[-,draw=green,thick] (z) -- (rhs);
      \draw[-,draw=purple,thick] (w2) -- (rhs);
      \draw[-,draw=cyan,thick] (x2) -- (rhs);
      \draw[-,draw=blue!50,thick] (y2) -- (rhs);
      \draw[-,draw=blue!50,thick] (z2) -- (rhs);
      \draw[-,thick] (c) -- (rhs);
      \draw[-,thick] (w) -- (w2);
      \draw[-,thick] (x) -- (x2);
      \draw[-,thick] (y) -- (y2);
      \draw[-,thick] (z) -- (z2);
    \end{tikzpicture}
    \caption{signed permutation symmetries}
    \label{fig:detectsignpermsym}
  \end{subfigure}
  \caption{Illustration of examplary symmetry detection graphs.}
  \label{fig:detectsym}
\end{figure}
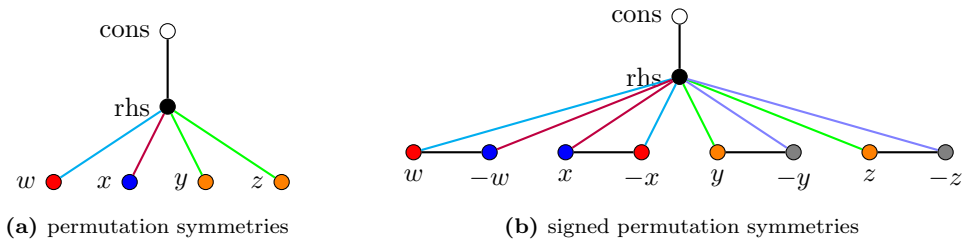

A possible symmetry detection graph for our exemplary \MILP is shown in
Figure~\ref{fig:detectpermsym}.
The nodes for variables~$y$ and~$z$ receive the same color since both have the
same objective coefficient and bounds; the remaining variable nodes receive
different colors, since they are not symmetric to each other.
Moreover, we introduce one node for the right-hand side, which is colored
according to the right-hand side coefficient.
The edges connect the variable nodes with the right-hand side node; they
are colored according to the coefficient of the corresponding variable in
the linear constraint.
This construction can easily be extended to general linear constraints,
see~\cite{Salvagnin2005}.

The symmetry detection graphs for individual constraints are then combined
into a single symmetry detection graph.
The callback \texttt{CONSGETPERMSYMGRAPH} provides a pointer to this graph
and different functions can be used add nodes and edges to this ``global''
symmetry detection graph.
To avoid re-defining variable nodes for different constraints, these nodes
cannot be added within the callback.
Instead, these nodes are defined centrally by \scip and are colored
according to the variable's type, its objective coefficient, and lower and
upper bound.
To make sure that only constraints of the same type can be symmetric to
each other (compare the permutation~$\pi$ above), every constraint should
add a ``constraint'' node to its local symmetry detection graph, which
serves as an identifier of the type of constraint (e.g., ``linear'',
``knapsack'', or ``SOS1'').

The functions for adding nodes to the global graph are
\texttt{SCIPaddSymgraphValnode()}, \texttt{SCIPaddSymgraphOpnode()}, and
\texttt{SCIPaddSymgraphConsnode()}.
The first function adds nodes that hold a numerical value, e.g., the
right-hand side node of a linear constraint.
The second function can be used to add ``operator'' nodes that allow to
model special relations between other nodes.
For example, in a nonlinear constraint such an operator could model
(nonlinear) functions that are applied to other nodes (e.g., variables).
Operators are encoded as integer numbers, i.e., the implementation needs to
make sure that only equivalent operators are assigned the same integer
value.
The third function adds a node that stores a pointer to the corresponding
constraint.

Edges can be added by \texttt{SCIPaddSymgraphEdge()}.
To add edges to variable nodes, the function
\texttt{SCIPgetSymgraphVarnodeidx()} can be used to get the index of the
variable node in the symmetry detection graph.

\paragraph{Detection of Signed Permutation Symmetries}
Let~$e^1,\dots,e^n$ be the standard unit vectors in~$\R^n$.
A signed permutation is a bijective map~$\gamma\colon\{\pm e^1,\dots,\pm e^n\} \to
\{\pm e^1,\dots,\pm e^n\}$ that satisfies~$\gamma(-e^i) = -\gamma(e^i)$ for
all~$i \in [n]$.
A signed permutation~$\gamma$ acts on a vector~$x \in \R^n$ as~$\gamma(x) =
(\sign(\gamma^{-1}(1))\, x_{|\gamma^{-1}(1)|}, \dots,
\sign(\gamma^{-1}(n))\, x_{|\gamma^{-1}(n)|})$, where~$\sign(\cdot)$ is the
sign function.
That is, it permutes the entries of the vector, but
it can also change the sign of some entries.
Signed permutation symmetries of an \MILP or \MINLP can be defined analogously
to permutation symmetries.
Such symmetries arise, e.g., in geometric problems like packing circles
into a box~\cite{CostaHansenLiberti2013}, where they model reflections along standard
hyperplanes.

To detect signed permutation symmetries, constraint handlers need to
implement the callback \texttt{CONSGETSIGNEDPERMSYMGRAPH}.
The functionality is analogous to the one of \texttt{CONSGETPERMSYMGRAPH}, however,
the automorphisms of the symmetry detection graph need to encode signed
permutation symmetries now.
This can be achieved by introducing not only a node for every variable~$x$,
but also for every negated variable~$-x$.
The colors of the negated variables are derived based on the negated
objective coefficient and negated bounds of variable~$x$.
Moreover, to indicate that~$x$ and~$-x$ form a pair of negated variables,
both must be connected by an edge in the symmetry detection graph.
Finally, the definition of signed permutations above requires that only
reflections along standard hyperplanes are allowed.
In \scip's implementation, however, also reflections along translated
standard hyperplanes can be detected.
This is achieved by not defining the colors based on the original variable
bounds, but for variables whose domain is translated to be centered at the
origin (except for semi-unbounded variables).
Consequently, also binary variables may admit signed permutation symmetries since
the variable domain is translated to~$\{-\frac{1}{2}, +\frac{1}{2}\}$.

Figure~\ref{fig:detectsignpermsym} shows the symmetry detection graph for our illustrative
example.
Next to the classical permutation symmetry that exchanges~$y$ and~$z$, also
the signed permutation that maps variable~$w$ onto
variable~$-x$ can be detected.

To create the symmetry detection graph for signed permutation symmetries
all functions to add nodes and edges as described above can be used.
To access the index of a negated variable,
\texttt{SCIPgetSymgraphNegatedVarnodeidx()} needs to be used.

Via the parameter \texttt{propagating/symmetry/symtype}
a user can select which type of symmetries is detected.
A value of~0 corresponds to permutation symmetries, and a value of~1 to
signed permutation symmetries.
By default, signed permutation symmetries are not detected, because
currently only basic symmetry handling techniques for such symmetries are
implemented.

\subsubsection{Symmetry Interfaces}
\label{sec:symmetryinterfaces}
To detect automorphisms of the previously mentioned symmetry detection
graphs, the previous version of \scip made use of the graph automorphism
software \bliss~\cite{bliss}, which is also shipped together with \scip.
\scipv also features interfaces to \solver{nauty/traces}~\cite{Nauty}.
Depending on which software package shall be used for symmetry detection, the \scip
make command takes \texttt{SYM=\{none,bliss,nauty\}} as argument.
Moreover, \scipv allows to make use of \sassy~\cite{Sassy},
which preprocesses the symmetry detection graphs to accelerate the
computation of symmetries.
\sassy can be used by compiling \scip with option
\texttt{SYM=\{sbliss,snauty\}}.

\subsection{Nonlinear Handlers}
\label{subsect:nlhdlrs}

Nonlinear constraints (see \eqref{eq:generalminlp}) are handled by the constraint handler \texttt{nonlinear} in \scip.
This constraint handler can delegate tasks on detecting and exploiting structure in algebraic expressions to specialized \textit{nonlinear handlers}, see~\cite{SCIP8} for details.
For example, the nonlinear handler for quotient expressions identifies expression of the form $vw^{-1}$ (where $v$ and $w$ can be arbitrary expressions) in nonlinear handlers and provides specialized bound tightening and linear under/overestimators for the function $(v,w)\mapsto vw^{-1}$.

With \scipv, a new nonlinear handler for signomial functions has been added and the nonlinear handler for quadratic expressions has been improved.
In addition, a new nonlinear handler callback has been added to request the linearization of an expression in a given solution point.
The nonlinear constraint handler can use this callback to tighten the linear relaxation when a new feasible solution has been found (parameter \texttt{constraints/nonlinear/linearizeheursol}, currently disabled by default).

\subsubsection{Signomial Handler}
\label{subsubsect:signomial}

An $n$-variate \emph{signomial term} is defined as $ x^{\alpha} =\prod_{j \in [n]} x_j^{\alpha_j}$, where  $\alpha \in \R^n$ and $x > 0$.  In general, the signomial term is nonconvex. In \scipv, a new nonlinear handler is implemented that generates cutting planes for signomial constraints.

Given $x^{\alpha}$, the handler aims at approximating the lifted set $S\defi \{(x,t) \in \R^n \times \R: t = x^{\alpha}  \}$, which is in general given by a constraint of the extended formulation (see the 8.0.0 release report~\cite{SCIP8} for details on extended formulations). It is easy to show that $S$ can be rewritten in the form
\begin{equation}
\label{eq.form1}
	S =  \{(u, v) \in \R_{+}^{h + \ell}:\, u^{\bar{\beta}} = v^{\bar{\gamma}}\}.
\end{equation}
with $\bar{\beta}, \bar{\gamma}$ containing only nonnegative entries.

Given a point $(\tilde{u},\tilde{v})$, the handler outer approximates either $S_1 \defi  \{(u, v) \in \R_{+}^{h + \ell}:\, u^{\bar{\beta}} \le v^{\bar{\gamma}}\}$ or $S_2\defi \{(u, v) \in \R_{+}^{h + \ell}:\, u^{\bar{\beta}} \ge v^{\bar{\gamma}}\}$ by checking which of the two sets does not contain $(\tilde{u},\tilde{v})$. More precisely, if $(\tilde{u},\tilde{v}) \notin S_1$, then the handler separates a linear valid inequality for $S_1$ that (possibly) cuts off $(\tilde{u},\tilde{v})$ and overestimates the signomial term; if $(\tilde{u},\tilde{v}) \notin S_2$, then the handler separates a linear valid inequality for $S_2$ that (possibly) cuts off $(\tilde{u},\tilde{v})$ and underestimates the signomial term.
In the case of an inequality constraint, only one of the sets $S_1$, $S_2$ correctly describes the feasible set of the constraint, and thus only one set is considered for separation.

The above formulation  exhibits symmetry between $u$ and $v$.  We only illustrate the method to approximate $S_1$, and the similar result applies to $S_2$ as well.

Since the signomial terms $u^{\bar{\beta}}, v^{\bar{\gamma}}$
are nonnegative over $\R_{+}^h$, $\R_{+}^\ell$, we can take any positive power $\eta \in \R_{>0}$ on both sides of \eqref{eq.form1} to obtain
\begin{equation}
\label{eq.sig_form2}
	S_1 =  \{(u, v) \in \R_{+}^{h + \ell}:\, u^\beta \le v^\gamma \},
\end{equation}
where $\beta \defi  \eta {\bar{\beta}}$, $\bar{\gamma} \defi \eta \bar{\gamma}$, and $\eta = 1 / \max(\sum_{j \in [h]} \abs{\bar{\beta}_j},  \sum_{j \in [\ell]} \abs{\bar{\gamma}_j})$. Thus, we have that $\max(\sum_{j \in [h]} \abs{\beta_j},  \sum_{j \in [\ell]} \abs{\gamma_j}) = 1$.

Moreover, we assume that the range  of $u$ is a hyperrectangle $U \subseteq \R^h_{>0}$, which is usually available as variable bounds in \scip.
The reformulated set enjoys two useful properties \cite{xu2022cutting}:  the terms $u^\beta, v^\gamma$ in \eqref{eq.sig_form2} are concave functions, and the convex envelope of $u^\beta$ over $U$ is vertex polyhedral. The signomial handler aims at (1) linearizing the convex envelope of $u^\beta$, i.e., an affine underestimator for $u^\beta$, and (2) linearizing the concave function  $v^\gamma$.

Let $Q$ be the vertices of $U$.  Since the convex envelope of $u^\beta$ over $U$ is vertex polyhedral, the separation of an affine underestimator $a \cdot u + b$ could be solved by an LP:
\begin{equation}
\label{eq.sig_sep}
	\max_{a \in \R^{h}, b \in \R} \{a \cdot \tilde{u}+b \st \forall q \in Q \; a \cdot q+b \le q^\gamma\},
\end{equation}
When $h = 1, 2$, the handler directly uses a closed form expression of the optimal $a$, $b$ without solving the LP \cite{xu2022cutting}.

Now, denote $g(v) := v^\gamma$. Since $g$ is concave, a straightforward way to overestimate it is to use the gradient at $\tilde{v}$ and obtain $g(\tilde{v}) + \nabla g(\tilde{v}) \cdot (v - \tilde{v}) $. Then, the separated valid inequality is of the form $a \cdot u + b \le g(\tilde{v}) + \nabla g(\tilde{v}) \cdot (v - \tilde{v})$. In the implementation, such inequality is further transformed into an overestimator of $x^\alpha$ through scaling.

Table \ref{tab:sigcut} shows the impact
of the signomial handler on \scip performance on 152 MINLPLib instances
that contain signomial terms. We report the ratio of shifted geometric means of time and nodes and the number of solved instances.
\begin{table}[h]
	\caption{Performance statistics of the signomial handler over \scip default.}
	\label{tab:sigcut}
	\centering
	\begin{tabular}{ccc}
		\toprule
		\multicolumn{3}{c}{152 selected MINLPLib instances} \\
		\midrule
		Time (s) & Nodes & Solved  \\
		\midrule
	\textbf{0.92}	 & \textbf{0.93} & \textbf{78 vs 75}  \\
		\bottomrule
	\end{tabular}
\end{table}
 
The computational cost of estimating a signomial term primarily stems from solving the  LP \eqref{eq.sig_sep}, whose size is exponential in $h$. An advanced parameter that governs the maximum allowable value of $h$ that the signomial handler can manage is \texttt{nlhdlr/signomial/advanced/maxnundervars}.
Currently, the signomial handler is disabled by default. Users can enable it via the \texttt{nlhdlr/signomial/enabled} parameter.

The current implementation of the signomial handler lacks a specialized bound tightening method for variables within signomial terms.
Given the critical impact of variable range ($U$) on cut quality, a further development of the handler involves refining these bounds through propagation techniques.

\subsubsection{Strengthening Cuts for Quadratic Expressions}
\label{subsubsect:quadratic}
To separate nonconvex quadratic constraints, the constraint handler \texttt{nlhdlr\_quadratic}, which was introduced in \scip 8.0, can generate intersection cuts by setting the parameter \texttt{nlhdlr/quadratic/useintersectioncuts = TRUE} (currently disabled by default).
Until now, the intersection cuts were built by using only the current \LP relaxation and the violated quadratic constraint.
To additionally leverage integrality information, \scipv allows to strengthen the cutting planes using \emph{monoidal strengthening} \cite{BALAS1980224} if some of the variables in the problem need to be integer.
As \cite{chmiela2023monoidal} showed, the strengthened intersection cuts significantly outperform the pure intersection cuts whenever monoidal strengthening can be applied.

\subsection{Primal Heuristics}
\label{subsect:heuristics}

\subsubsection{Indicatordiving}
\label{subsubsect:indicatordiving}
\emph{Semi-continuous} variables are variables that take either the value 0 or any value within a specific range:
\begin{equation*}
	x \in \{0\} \cup [\lb,\ub] \text{ with } 0 < \lb \leq \ub \text{ and } \ub \in \R_{+}\cup\infty.
\end{equation*}
Such variables are used, for example, in modeling supply chains where a facility either can produce nothing or, if enabled, has to produce at least an amount $\lb$.

Semi-continuous variables can be formulated in \scip with an additional binary variable $z\in\{0,1\}$ as
\begin{align}
	x \in [0,\ub], \nonumber\\
	\lb z \leq x,\label{eq:indicatordiving:reform:lin}\\
	z = 0 \implies x \leq 0.\label{eq:indicatordiving:reform:ind}
\end{align}
Thereby,~\eqref{eq:indicatordiving:reform:ind} is a so-called indicator constraint, that is, $x\leq0$ must hold if $z=0$.
Linear constraint~\eqref{eq:indicatordiving:reform:lin} models the lower bound on $x$.

If $\ub$ is finite, one can reformulate the indicator constraint~\eqref{eq:indicatordiving:reform:ind} with a linear big-$M$ constraint, such as $x \leq \ub z$.
If $\ub$ is infinite, this is not directly possible. However, one could add an artificial upper bound with the risk of cutting off optimal solutions and causing numerical issues due to the large upper bound $M$.
A new diving heuristic, \texttt{indicatordiving}, has been developed to find
solutions also in the presence of indicator constraints modeling semi-continuous variables with infinite upper bounds~$\ub$.

Diving heuristics iteratively fix variables and solve the modified \LPs simulating a depth-first-search in an auxiliary tree.
A description of the generic diving procedure used in \scip can be found in~\cite{SCIP6}.
Other diving heuristics in \scip typically take only integer variables with fractional \LP solution value into account.
In contrast, \texttt{indicatordiving} additionally examines all binary indicator variables $z$ corresponding to violated indicator constraints and which are integral in the \LP solution but not fixed already.

Each such variable is assigned a score to determine the variable that should be fixed next.
For indicator variables, the score is given by
\begin{equation*}
\phi \coloneqq
	\begin{cases}
		-1, &\text{ if } \hat{x} \in \{0\} \cup [\lb,\ub],\\
		100 \cdot (\lb - \hat{x})/\lb, &\text{ if } \hat{x} \in (0,\ub),
	\end{cases}
\end{equation*}
where $\hat{x}$ is the current \LP solution value.
The indicator variable $z$ with the highest score gets fixed to $1$ if the \LP value $\hat{x}$ is at least 50\% of the lower bound $\lb$.
Otherwise, it gets fixed to $0$.
As soon as all the indicator variables are integral in the \LP solution or all the indicator constraints are fulfilled,
other candidate variables are considered, for which the score and rounding direction of \texttt{farkasdiving} are used.

\subsubsection{Extension of Dynamic Partition Search}
\label{subsubsect:dps}
Since \scip~7.0, the decomposition information can be passed to \scip in addition to the instance, which can be leveraged in heuristics, for example.
A detailed description of decompositions and their handling in \scip can be
found in the release report for that version~\cite{SCIP7}.

Dynamic Partition Search (\texttt{DPS}) introduced in \scip~8.0 is a heuristic that requires a decomposition.
\texttt{DPS} splits an \MILP~\eqref{eq:generalmilp} into several subproblems according to a decomposition.
Thereby, the linking constraints and their right-/left-hand sides are also split by introducing new parameters $p_q$ for each block $q$, called \emph{partition}.
Such a partition is central to the \texttt{DPS}.
When the heuristic is called during node preprocessing, the partition is initialized with a uniform distribution of the constraint sides over the blocks.

In \scipv, the \texttt{DPS} has been extended with an option to get called at the end of the node processing and, therefore, can use the \LP solution for initialization of the partition.
The parameter \texttt{heuristics/dps/timing} controls the calling point and, thus, the initialization.
A detailed description of \texttt{DPS} and additional heuristics exploiting decomposition information can be found in~\cite{Halbig_DecHeur_2023}.

\subsubsection{Learning to Control Primal Heuristics Online}
\label{subsubsect:adaptiveheuristics}
Since the performance of heuristics is highly problem-dependent and most of them can be very costly, it is necessary to handle them strategically.
Thus, sometimes it is preferable to have dynamic, self-improving procedures rather than relying on static methods to control primal heuristics.

\scip has already used two adaptive heuristics in previous versions that use bandit algorithms to decide which heuristics to additionally run:
Adaptive Large Neighborhood Search (ALNS)~\cite{hendel2018alns} and Adaptive Diving~\cite{hendel2018bandit}.
Building upon this, \scipv now includes a general online learning approach~\cite{chmiela2023heuristic}, which dynamically adapts the application of primal heuristics to the unique characteristics of the current instance.
In particular, both \textit{Large Neighborhood Search} and \textit{Diving} heuristic types are controlled together, making this the first work where two different classes of heuristics are treated simultaneously by a single learning agent.
In \scipv, this is implemented as a heuristic called \texttt{scheduler}.
Since this framework was designed to replace the classical heuristic handling as a whole as opposed to being run as an additional heuristic, \texttt{scheduler} is disabled by default in \scipv.

\subsection{Cutting Planes}
\label{subsect:cuts}
This section discusses the updates to the separation routine in \scip, both for cut generation and cut selection.
Separation in \scip is performed in rounds.
In a round, various valid inequalities that cut off the current LP relaxation's fractional solution are generated and stored either in a global cut pool (for cuts that are valid globally) or a separation store (for cuts that are valid only locally).
Then, these cuts are filtered and added to the LP relaxation before re-solving the relaxation and proceeding to the next round of separation.
Note that other components of \scip such as branching maybe executed in between two separation rounds.
These rounds are performed until a stopping criterion is met (e.g., maximum number of rounds or cuts added, or dual bound of the relaxation stalling).
\subsubsection{Cut Generation}
\label{subsubsect:cutgeneration}
\paragraph{Lagromory Separator}
There are two potential issues with the round-based approach mentioned above.
\begin{enumerate}
 \item The generation of higher-ranked cuts, e.g., higher-ranked Gomory Mixed-Integer (GMI) cuts~\cite{cornuejols2013safety}, may result in numerical troubles during the solving process.
 
 \item Multiple rounds of separation may be needed to achieve dual bound improvement in the presence of dual degeneracy~\cite{gamrath2020exploratory}.
\end{enumerate}
The first issue is addressed in the literature via a relax-and-cut framework-based separation techniques~\cite{guignard1998efficient, lucena2005non, cavalcante2008relax, fischetti2011relax}.
The second issue has received less attention until now despite being critical to the usefulness of separation in the solvers.
The new separator in \scip, the Lagromory separator, addresses both these issues.
It is a relax-and-cut framework-based separator which is built based on the separation technique presented in~\cite{fischetti2011relax}.

In the basic version of the relax-and-cut framework, which is also discussed in~\cite{fischetti2011relax}, when a separator is called at a node with fractional \LP solution, certain cuts are generated but not added directly to the \LP relaxation.
These cuts are added to the objective function of the node \LP in a Lagrangian fashion using Lagrangian (penalty) multipliers.
Then, this Lagrangian dual problem is solved via an iterative approach by updating the Lagrangian multipliers in every iteration, requiring an \LP solving in every iteration.
When an \LP is solved in an iteration, it is equivalent to exploring a new basis of the node \LP.
Then, additional cuts are generated with respect to this newly explored basis and are added to the objective function of the node \LP again.
This procedure is repeated until certain termination criterion is met.
While this approach in~\cite{fischetti2011relax} was proven to improve the dual bound at a given single node, it turned out to be ineffective in the context of the entire branch-and-cut tree.
To overcome this crucial hurdle, the Lagromory separator in \scip also implements various novel enhancements.
\begin{itemize}
  \item Theoretical enhancements include stabilization and regularization of the vector of Lagrangian multipliers.
  This vector is integral to and iteratively updated in the relax-and-cut framework.
  Stabilization using a core vector is an essential component in the literature of decomposition methods such as the Benders' and Dantzig-Wolfe decompositions.
  Regularization of vectors (e.g., by projecting the vectors into $\ell_1$, $\ell_2$, etc, norm balls) is a commonly applied technique in the literature of nonlinear optimization.
  \item Computational enhancements include the threshold for dual degeneracy beyond which the separator is executed; the working limits on the number of \LP iterations, number of cuts generated per explored basis of the \LP relaxation; etc.
\end{itemize}

The separator was tested on the MIPLIB 2017 benchmark library~\cite{miplib2017}.
It speeds up the solving process of harder instances that require at least 1000 seconds for solving to optimality.
On the other hand, it increases the solving time for many easier instances resulting in no overall improvement of the default \scip performance.
The Lagromory separator is OFF by default due to this reason.
An interested user may switch it ON by changing the parameter \texttt{separation/lagromory/freq} to a non-negative number.

\subsubsection{Cut Selection}
\label{subsubsect:cutselection}
The cut selector plugin introduced in \scip~8.0 enabled multiple research directions on the problem of cut selection, which recently underwent a revived scrutiny.
The previously hard-coded algorithm was replaced by the default cut selector \texttt{cutsel/hybrid}, which in particular scores cuts with a weighted sum of four criteria: efficacy, integer support, objective parallelism, and directed cutoff distance (with the last one having a zero default weight).
The scored cuts are then filtered iteratively by an orthogonality criterion.
This removes all non-orthogonal cuts (within some tolerance), starting from the highest-scoring cut until every cut has been processed, and the remaining set of cuts is pairwise near-orthogonal.
The importance of cut selection and some limitations of the current criteria were in particular shown theoretically and computationally in \cite{OJMO_2023__4__A5_0,turner2023cutting}.

These lines of work point to the conclusion that the current cut selection algorithm fails to capture and adapt to important instance properties.
In \cite{turner2023context}, a new cut selector coined \emph{ensemble} was developed to capture more of these properties than the current weighted sum of the default \texttt{cutsel/hybrid}.
The three core aspects of the cut selection loop described above are the \emph{filtering} of cuts, their \emph{scoring}, and the \emph{stopping criteria} for the cut loop.
All three aspects were extended and modified in the new ensemble cut selector, resulting in a large number of parameters that were activated and adjusted through the blackbox hyperparameter optimization tool SMAC \cite{JMLR:v23:21-0888}.
These parameters include pseudo-cost based cut scoring, sparsity based cut scoring, density based cut filtering, parallelism based scoring penalties, and a stopping criterion based on the number of nonzeros in the added cuts.
Because of the filter on cut density, \texttt{cutsel/ensemble} performs unreliably on MINLPs, where some nonlinear constraints demand the use of high-density cuts for good solver performance.
It is therefore applied with a lower priority than \texttt{cutsel/hybrid} which remains the default, but can be activated by changing \texttt{cutselection/ensemble/priority} to a number greater than 8000, the value of \texttt{hybrid}.

\begin{figure}
    \centering
    \begin{minipage}{.48\textwidth}
    \centering
    \includegraphics[width=\textwidth, keepaspectratio=true]{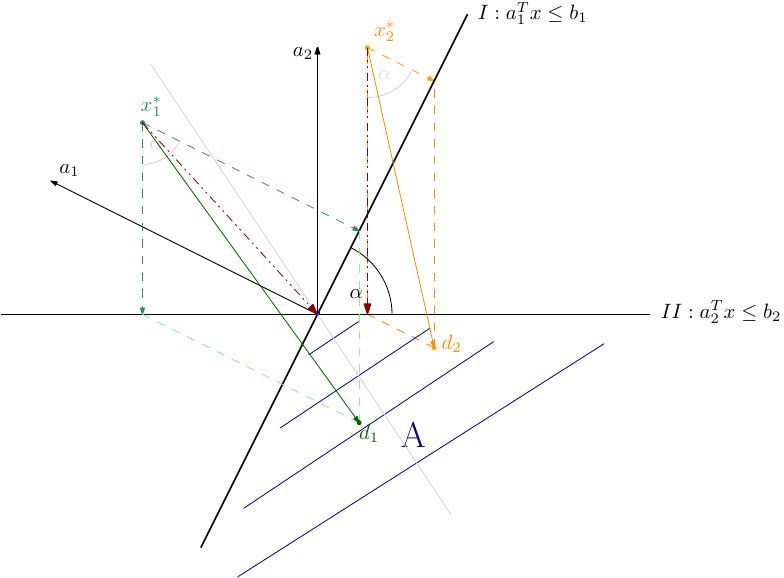}
    (a)
    \end{minipage}%
    \begin{minipage}{.48\textwidth}
    \centering
    \includegraphics[width=\textwidth, keepaspectratio=true]{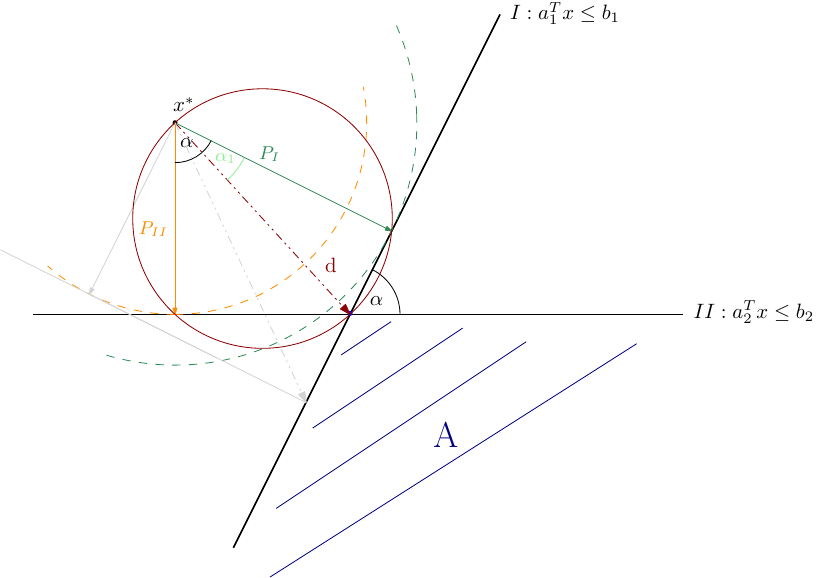}
    (b)
    \end{minipage}
\captionof{figure}{(a) Pairwise cut efficacies, shown for two cases of \LP solutions in yellow and green for cutting planes I and II. (b) Orthogonality-based cutfiltering, shown by the offset between the pairwise efficacies in red and gray of cut I paired with either cut II or the orthogonally rotated cut II.}
\label{fig:near_ortho_filter}
\end{figure}

Further, a novel dynamic filtering method is introduced in \texttt{cutsel/dynamic}, which aims to enhance the current near-orthogonal threshold methodology used in the \texttt{cutsel/hybrid} selector.
The geometric rationale for this approach can be understood by considering how the \LP improvements for a pair of cutting planes, in terms of efficacy, are influenced by the position of their intersection relative to the current \LP optimal solution $x^*$ (refer to Figure \ref{fig:near_ortho_filter}a).

This is exemplified by demonstrating that a mere summation
of individual cut efficacies can be misleading in evaluating the actual efficacy of a pair of cuts.
The true efficacy, indicated by red arrows in Figure \ref{fig:near_ortho_filter}a, represents the minimal distance between the feasible region $A$ and the current \LP solution post-application.

Figure \ref{fig:near_ortho_filter}a illustrates a fundamental issue when using pure efficacy to evaluate multiple cuts simultaneously. 
This underscores the motivation for orthogonality-based filtering as showcased in Figure~\ref{fig:near_ortho_filter}b.
Specifically, if the current \LP solution is outside the fan formed by the intersection of the cuts (yellow), then one of the cuts becomes entirely ineffective, provided that the vertex created by this pair is not the optimal point in the subsequent iteration.
Conversely, if the \LP solution is within this fan, the aggregated efficacy of the cuts does not accurately reflect their true effectiveness.
This discrepancy is due to the degree of non-orthogonality between the cuts.

To address these limitations, the stringent orthogonality threshold of the default selector was relaxed. 
The new dynamic criterion is depending on individual cut efficacies and aims to position the \LP solution within the intersection fan of the cut pairs (as depicted in Figure \ref{fig:near_ortho_filter}b).
Additionally, users can specify a minimum efficacy improvement relative to the previous cutting plane via the \textit{mingain} parameter.
The \textit{filtermode} parameter 'f' in \texttt{cutsel/dynamic} facilitates rescoring of cuts between filtering rounds, using the pairwise efficacy instead of the usual scoring mechanism.
Preliminary test results have not indicated a general improvment over the default cut selection; hence, this selector is assigned a lower priority than \texttt{cutsel/hybrid} as well.

\subsection{Branching}
\label{subsect:branching}

\scip9 introduces a new branching criterion both implemented as a stand-alone rule and integrated within the default hybrid branching rule \cite{achterberg2009hybrid}.
For a more thorough overview of the results and idea presented in this section, see \cite{turner2023branching}.

\subsubsection{GMI Branching}
\label{subsubsec:GMI}

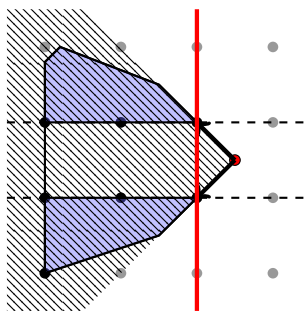
\begin{figure}[h]
	\centering
	\begin{tikzpicture}[scale=1]
		\centering
		\foreach \x/\y in {0/0, 0/1, 0/2, 0/3, 1/0, 1/1, 1/2, 1/3, 2/0, 2/1, 2/2, 2/3, 3/0, 3/1, 3/2, 3/3}
		{
			\fill[opacity=0.4] (\x,\y) circle (2pt);
		}
		\foreach \x/\y in {0/0, 0/1, 0/2, 1/1, 1/2, 2/1, 2/2}
		{
			\fill[opacity=1] (\x,\y) circle (2pt);
		}
		\node (a) at (2.5,1.5) {};
		\node (b) at (1.5,0.5) {};
		\node (c) at (0,0) {};
		\node (d) at (0,2) {};
		\node (e) at (0,2.8) {};
		\node (f) at (0.2,3) {};
		\node (g) at (1.5,2.5) {};
		
		\node (d11) at (2, 1) {};
		\node (d12) at (0, 1) {};
		\node (d21) at (2, 2) {};
		\node (d22) at (0, 2) {};
		
		\node (c1) at (0.5,-0.5) {};
		\node (c2) at (-0.5,-0.5) {};
		\node (c3) at (-0.5,3.5) {};
		\node (c4) at (0.5,3.5) {};
		
		\draw[thick] (a.center) -- (b.center) -- (c.center) -- (d.center)  -- (e.center) -- (f.center) -- (g.center) -- cycle;
		\filldraw[thick,fill=blue!60,fill opacity=0.4] (d11.center) -- (b.center) -- (c.center) -- (d12.center) -- cycle;
		\filldraw[thick,fill=blue!60,fill opacity=0.4] (d21.center) -- (d22.center) -- (d.center) -- (e.center) -- (f.center) -- (g.center) -- cycle;
		\fill[draw=black,fill=red] (a.center) circle (2pt);
		
		\draw[ultra thick,->] (a.center) -- (2,2);
		\draw[ultra thick,->] (a.center) -- (2,1);
		
		\fill[pattern = north west lines] (a.center) -- (c1.center) -- (c2.center) -- (c3.center) -- (c4.center) -- cycle;
		
		\draw[thick,dashed] (-0.5,1) -- (3.5,1);
		\draw[thick,dashed] (-0.5,2) -- (3.5,2);
		\node (f) at (1.2,1.1) {};
		
		\node (g) at (0.8,0.2) {};
		\node (h) at (0.8,2) {};
		
		\draw[ultra thick, red] (2, -0.5) -- (2, 3.5);
		\node (xxx) at (2,3.5) {};
		\node (xxy) at (2,-0.5) {};
	\end{tikzpicture}
\caption{Example intersection cut that is also a split cut.}
\label{fig:split_cut}
\end{figure}

GMI cuts~\cite{gomory1960algorithm} are a standard tool for solving MILPs.
It is recommended to see \cite{andersen2005reduce,cornuejols2008valid,turner2023branching} for an overview of how GMI cuts are derived and applied in practice.
The motivation to use GMI cuts to influence branching decisions stems from the observation that GMI cuts are derived from a split disjunction. A \emph{split disjunction} is defined by an integer $\pi_0 \in \Z$ and an integral vector $\pi \in \Z^{p} \times \{0\}^{n - p}$, which has zero entries for coefficients of continuous variables.
The disjunction is then given by:
\[
	\{ x \in \R^n \mid \pi_0^T x \leq \pi_0 \} \cup \{ x \in \R^n \mid \pi_0^T x \geq \pi_0 + 1 \}.
\]
Since no integer point is contained in the split $\{ x \in \R^n \mid \pi_0 < \pi_0^T x < \pi_0 \}$, all feasible points lie in the above split disjunction.
Thus, a cut that only cuts off points from the continuous relaxation that are inside the split is valid for the original problem. Such cuts are called \emph{split cuts} and GMI cuts are a special case of this family of cutting planes.
An example split with a valid split cut is visualized in Figure~\ref{fig:split_cut}.

To best link GMI cuts and branching, note that for $\pi = e_i$ and $\pi_0 = \lfloor \bar{x}_i \rfloor$, where $e_i$ denotes the $i$-th unit vector and $\bar{x}$ the current \LP solution, the split corresponds exactly to the region that is excluded by branching on a fractional variable $x_i$, $i \in \mathcal{I}$. The simple logic of this branching rule is then the following: Variables whose splits would generate a deep cut might also be good branching candidates. Therefore, the new Gomory branching rule generates all GMI cuts (or up to some maximum number of candidates if set), and branches on the variable whose associated split produces the most efficacious GMI cut.

\subsubsection{Using GMI Cuts for Reliability Pseudo-Cost (Hybrid) Branching}
\label{subsubsec:GMI_improve}

The default branching rule of \scip \cite{achterberg2009hybrid} has been extended to include two new terms in the weighted sum scoring rule. Currently, the variable that is branched on by \scip is the one with the highest branching score (ignoring cases of epsilon-close results).
The score for each variable is computed via a weighted sum rule that combines the following measures:
\begin{itemize}
	\item The frequency that the variable appears in a conflict
	\item The average length of the conflicts that the variable appears in
	\item The frequency that the variable's branching history has resulted in other variables becoming fixed
	\item The frequency that the variable's branching history has resulted in infeasible subproblems
	\item The number of NL constraints the variable features in
	\item The pseudo-cost associated with the variable from its history of previous branching decisions.
\end{itemize}

This weighted sum rule is appended with two terms controlled by the weight parameters \texttt{gmiavgeffweight} and \texttt{gmilasteffweight}.
In a separation round, \scip now stores the normalized efficacy of generated GMI cuts.
The normalization is a simple division by the largest efficacy of any GMI cut generated in the separation round,
resulting in a value in the range $[0,1]$ for each variable from which the tableau row produced a GMI cut.
For each variable in the problem, \scip now records a running average of the normalized GMI cut efficacies from tableau rows associated with the respective variable.
Additionally, \scip now also records the last normalized value for each variable. 
The parameter \texttt{gmiavgeffweight} is then the weight of the running average in the weighted sum rule for each variable. The parameter \texttt{gmiavgeffweight} is equivalently the weight of the most recent recorded value for each variable.
By default, \texttt{gmiavgeffweight} is set to 0, and \texttt{gmilasteffweight} is set to $10^{-5}$, effectively acting as a tie-breaker.

\subsection{LP Interfaces}
\label{subsect:lpi}

\paragraph{\highs LPI}\label{subsect:highs}
\scipv provides the possibility of using the open-source \LP solver \highs\cite{HuangfuHall15} \footnote{Available in source code at \url{https://github.com/ERGO-Code/HiGHS}.}
.
The interface provides the basic functionality, yet it does not fully exploit all capabilities of \highs.

\subsection{Technical Improvements}
\label{subsect:further}

\paragraph{AMPL .nl reader}
Added support for logical constraints in binary variables and basic logical operators (and, or, not, equal).

\paragraph{OBBT propagator}
Variables of linear constraints that are controlled by indicator constraints can now also be taken into account for bound tightening.
This feature is disabled by default, but can be enabled via parameter \texttt{propagating/obbt/indicators}.

\section{SoPlex}
\label{sect:soplex}

Most importantly, \soplexv now supports \emph{incremental precision boosting}~\cite{Espinoza2006} for solving \LPs exactly over the rational numbers, in addition to and in combination with the existing \LP iterative refinement approach~\cite{GleixnerSteffyWolter2016,GleixnerSteffy2020}.
The algorithm for exact solving can be selected using the boolean parameters \texttt{precision\_boosting} and \texttt{iterative\_refinement}.
By default, both are set to \texttt{true}, in which case \soplex uses a combined algorithm with an outer precision boosting loop and an inner iterative refinement loop.
For further details and computational experiments we refer to~\cite{EiflerThouveninGleixner2023}.
The new default for exact, rational \LP solving increases the robustness of the algorithm on numerically challenging problems and allows to solve more problems exactly.
Furthermore, a new Python interface for \soplexv, called \pysoplex, has been developed, see~\cref{subsect:pysoplex}.

\section{\papilo}
\label{sect:papilo}

\papilo, a C++ library, is a solver-independent presolving library that provides presolving routines for \MIP and \LP and is part of the presolving routines in \scip.
It also supports multi-precision arithmetic, which makes \papilo an essential part of the presolving process of \exscip~\cite{exscipgithub}, the numerically exact version of \scip~\citep{CookKochSteffyWolter13Hybrid,EiferGleixner21ComputationalStatusUpdate,EiflerGleixner2023}.

\papilov now supports \emph{proof logging} as a new feature, i.e., the generation of machine-verifiable certificates by a solver in order to prove the correctness of its computation.
Proof logging was originally introduced by the SAT community to ensure the correctness of a solver's computation, since even state-of-the-art solvers falsely claim infeasibility or optimality or return infeasible ``solutions''~\cite{AkgunGentJeffersonMiguelNightingale18Metamorphic,CookKochSteffyWolter13Hybrid,Klotz2014,Steffy2011,GillardSchausDeville19SolverCheck}.
Examples of such proof formats are DRAT~\cite{WetzlerHeuleWarren14DRAT,HeuleWarren13Trimming,HeuleWarren13Verifying}, GRIT~\cite{CruzMarquesSchneider17EfficientCertified}, and LRAT~\cite{CruzHeuleKaufmannSchneider17EfficientCertified}.
\exscip has adapted proof logging to certify the branch-and-cut process using the \vipr format~\cite{CheungGleixnerSteffy17Verifying}, but since \vipr currently does not provide the necessary functionality to verify presolving reductions, \exscip only prints the certificate for the presolved problem.

As a first step, \papilov provides the ability to generate proofs for presolving of binary programs in the \veripb format\footnote{available at \url{https://gitlab.com/MIAOresearch/software/VeriPB}}, which was developed for \emph{Pseudo-Boolean} (PB) problems~\cite{BGMN23Dominance,GMNO22CertifiedCNFencodingPB}.
\veripb readily supports a \emph{reverse unit propagation} rule and a \emph{redundancy-based strengthening} rule for verifying dual arguments, and in our effort to certify presolving transformations, it has recently been extended by an \emph{objective update rule} to support modification of the objective during presolving.
We refer to \cite{HoenOertelGleixnerNordstromCertifying24} for a detailed explanation of how each presolving reduction in \papilo can be certified using the \veripb language.
In order to print a certificate in \papilo the boolean parameter \verb|verification_with_veripb| must be set to true.
Since \veripb only supports PB problems, proof logging is currently only supported for this problem class.

\Cref{tab::papilo::overhead_obju_diff} reports the performance impact of proof logging.
These experiments are based on the selection of binary programs from \miplib2017~\cite{miplib2017}, called \miplib01~\cite{MIPLIB01}, and the instances of the \nbsc{PB16} competition~\cite{PB16}, each split into optimization (opt) and decision (dec) instances.
We only exclude the large-scale instances \nbsc{ivu06-big} and \nbsc{supportcase11} with a runtime of more than 2 hours in \papilo.
Times are aggregated using the geometric mean shifted by 1~second.
The overhead of proof logging ranges from 27\% to 54\% on both test sets.
For 99\% of the decision instances, the \emph{overhead per applied transaction} is below $0.186$~milliseconds on both test sets. 
This shows the viability of proof logging in practice especially considering that proof logging runs sequentially while the presolvers in \papilo run in parallel.

\begin{table}[ht]
	\centering
	\small
	\begin{tabular*}{\textwidth}{@{}l@{\;\;\extracolsep{\fill}}lrrr}
		\toprule
		Test set & instances &  default (in seconds) & w/proof log (in seconds) & relative \\
		\midrule
		\nbsc{PB16}-dec & 1398 & 0.050 & 0.077 & 1.54\\[-0.5ex]
		\miplib01-dec&  295 & 0.498 & 0.631 &  1.27\\[-0.5ex] 
		\nbsc{PB16}-opt &  532 & 0.439 & 0.565 &  1.29 \\[-0.5ex]
		\miplib01-opt&  144 & 0.337 & 0.473 & 1.40 \\[-0.5ex] 
		\bottomrule
	\end{tabular*}
	\caption{Runtime comparison of \papilo with and without proof logging.}
	\label{tab::papilo::overhead_obju_diff}
\end{table}

\section{Interfaces}
\label{sect:interfaces}

\subsection{AMPL}
\label{subsect:ampl}
The AMPL\footnote{\url{https://www.ampl.com}} interface of \scip now supports parameters specified in AMPL command scripts via option \texttt{scip\_options}.
The value of \texttt{scip\_options} is expected to be a sequence of parameter names and values, separated by a space, for example
\begin{verbatim}
  option scip_options 'limits/time 10 display/verblevel 1';
\end{verbatim}

 \subsection{JSCIPOpt}
 \label{subsect:jscipopt}
The Java interface to \scip, \jscipopt, is again actively maintained and requires a C++ compiler (rather than just a C compiler) to compile. The following functionality of the {\scip} C API (or the \texttt{objscip} C++ API) is newly available from Java:
\begin{itemize}
 \item branch priorities,
 \item concurrent solving,
 \item changing the objective coefficient of a variable (contributed by the GitHub user xunzhang (Hong Wu)),
 \item message handlers (using the \texttt{objscip} API -- the C++ interface class can be transparently subclassed from Java),
 \item getting the current dual bound,
 \item getting the solving status (contributed by the GitHub user patrickguenther),
 \item interrupting the solving process,
 \item creating a partial solution (contributed by the GitHub user fuookami (Sakurakouji Sakuya)).
\end{itemize}
In addition, the following bugs were fixed:
\begin{itemize}
 \item running SWIG during the compilation only worked on *nix systems due to the unnecessary use of the external POSIX command \texttt{mv},
 \item building against {\scip} 8.0.0 or newer was failing due to a missing \texttt{\#include} statement,
 \item even when defining \texttt{SCIP\_DIR} at build time, a different {\scip} could be silently used if the passed \texttt{SCIP\_DIR} was not valid for some reason (this is now an error),
 \item a \jscipopt library dynamically linked to the \scip library was not binary-compatible with different versions of \scip due to the use of macros hardcoding structure layouts (macros now are avoided unless \scip is statically linked into \jscipopt),
 \item the \texttt{SCIP\_Longint} type was incorrectly mapped to a 32-bit Java type (now correctly mapped to a 64-bit Java type).
\end{itemize}

\subsection{PySCIPOpt}
\label{subsect:pyscipopt}
The Python interface to \scip, \pyscipopt \cite{MaherMiltenbergerPedrosoRehfeldtSchwarzSerrano2016}, is now automatically shipped with a standard installation of \scip when installed using PyPI. This automatic \scip installation is currently available only for machines running x86\_64 architecture. The Python versions and OS combinations supported include CPython 3.6+ for manylinux2014 (includes Ubuntu / Debian) and MacOS, and CPython 3.8+ for Windows. Linking \pyscipopt against a custom installation of \scip is still possible and encouraged, however now requires cloning \pyscipopt's repository, available on GitHub \cite{pyscipoptgithub}, and installing from source.

A new Python package, \pyscipoptml \cite{turner2023pyscipopt}, is now available. The package uses \pyscipopt to automatically formulate machine learning models into MIPs. This functionality allows users to easily optimize MIPs with embedded ML constraints, simplifying the process of deciding on a formulation and extracting the relevant information from the machine learning model.

\subsection{\scipjl}
\label{subsect:scipjl}
\scipjl is the Julia interface to \scip and provides access to the solver in two ways. First, it provides an access to all functions of the C interface mirrored by \solver{Clang.jl} and accessed via Julia \texttt{ccall}.
Second, it exposes a high-level interface implementing the \solver{MathOptInterface} API \cite{legat2022mathoptinterface} and callable, e.g., through the \solver{JuMP} modeling language \cite{lubin2023jump}.
The high-level interface now includes an access to the heuristic, branching, and cut selection plugins, making them available in an idiomatic Julia style, in addition to the constraint handler and separator plugins.
The heuristic and cut generation plugins are also available through the standardized \solver{MathOptInterface} callback mechanism.

\subsection{russcip}
\label{subsect:russcip}
With \scipv, we introduce the first version of the Rust interface for \scip, \russcip \cite{russcipgithub}.
The interface builds on Rust's solid foundation for type and memory safety. Being a system's programming language,
it allows for low-level access to the C-API of \scip, binding directly without the need for copying data across the language barrier. The interface is split into two parts: an unsafe part, which provides full access to the C-API through the module \texttt{ffi} and a limited but safe wrapper that allows access to part of the API. Currently, the following plugins are implemented on the safe interface that allow you to addi custom branching rules, primal heuristics, and variable pricers, and control \scip through event handlers, all designed to guarantee compliance with \scip's return types at compile-time.

The interface is still in its early stages and we are working on adding more plugins and improving the ergonomics of the safe interface. The following is a small example of how to use \russcip:

\begin{verbatim}
use russcip::prelude::*;

fn main() {
    // Create model
    let mut model = Model::default()
    .hide_output()
    .set_obj_sense(ObjSense::Maximize);

    // Add variables
    let x1 = model.add_var(0., f64::INFINITY, 3., "x1", VarType::Integer);
    let x2 = model.add_var(0., f64::INFINITY, 4., "x2", VarType::Integer);

    // Add constraint "c1": 2 x1 + x2 <= 100
    model.add_cons(vec![x1.clone(), x2.clone()],
     &[2., 1.], -f64::INFINITY, 100., "c1");

    let solved_model = model.solve();

    let status = solved_model.status();
    println!("Status:{:?}", status);

    let obj_val = solved_model.obj_val();
    println!("Objective:{}", obj_val);

    let sol = solved_model.best_sol().unwrap();
    let vars = solved_model.vars();

    for var in vars {
        println!("{}={}",&var.name(), sol.val(var));
    }
}
\end{verbatim}

Further examples related to defining and using custom plugins can be found in the repository's \cite{russcipgithub} tests.

\subsection{SCIP++}
\label{subsect:scippp}
\scippp is a C++ wrapper for \scip's C interface.
It automatically manages the memory, provides a simple interface to create linear expressions and inequalities, and provides type-safe methods to set parameters.
It can be used in combination with \scip's C interface, especially for features not yet present in \scippp.
The following is a small example.

\begin{verbatim}
#include <scippp/model.hpp>
using namespace scippp;
int main() {
    Model model("Simple");
    auto x1 = model.addVar("x_1", 1);
    auto x2 = model.addVar("x_2", 1);
    model.addConstr(3 * x1 + 2 * x2 <= 1, "capacity");
    model.setObjsense(Sense::MAXIMIZE);
    model.solve();
}
\end{verbatim}

\subsection{PySoPlex}
\label{subsect:pysoplex}
\pysoplex~\cite{pysoplexgithub} is a newly-introduced Python wrapper for \soplex's C interface.
The installation process is similar to that of \pyscipopt, so, a user needs to install \soplex first, set the \texttt{SOPLEX\_DIR} environment variable, and then install the \pysoplex wrapper.
The installation process has been successfully tested on Linux and Mac OS platforms.
The following is a small example.
\begin{verbatim}
import pytest
from pysoplex import Soplex, INTPARAM, BOOLPARAM, VERBOSITY

# create solver instance
s = Soplex()

# read instance file, solve LP, and get objective value
success = s.readInstanceFile("PATH_TO_INSTANCE.mps.gz")
# specify "lifting" parameter
s.setBoolParam(BOOLPARAM.LIFTING, 1)
# specify "verbosity" level
s.setIntParam(INTPARAM.VERBOSITY, VERBOSITY.ERROR)
s.optimize()
obj_val = s.getObjValueReal()
print(obj_val)    
\end{verbatim}

\section{The UG Framework}
\label{sect:ug}
\ug was originally designed to parallelize powerful state-of-the-art branch-and-bound based solvers (we call these ``\emph{base solvers}'').
Two of the most intensively developed parallel solvers are \fscip (for a shared memory computing environment) and
\pscip (for a distributed computing environment), both using \scip as the base solver.
\pscip solved two open instances on a supercomputer from \miplib (MIPLIB2003) for the first time in 2010~\cite{ParaSCIP-first}.
To achieve this, supercomputer jobs had to be restarted frequently from snapshots of the branch-and-bound tree.
To verify the results, we aimed to solve the instances with a single job on the supercomputer, which required the development of new features and intense debugging of \pscip.
Since debugging on distributed environments is inefficient, \fscip was developed, which has the same parallelization algorithms as that of \pscip (since \ug abstracts the parallelization library), but can run on a single PC.
The results of \fscip were first presented in the MIPLIB2010 paper~\cite{miplib2010} (therein, \fscip is referred to as UG[SCIP/SPX]).
Even though \fscip was already working in 2010, the \fscip paper~\cite{Shinano18fiber} was submitted only 3 years later, since the software went through intensive tests.
The supplement of the \fscip paper includes only a small fraction of the computational results we had conducted.
Due to the major debugging effort of \ug and \scip via \fscip, \pscip could solve more than 20 open instances~\cite{Shinano-ParaSCIP} from \miplib and none of these results have been proven wrong so far.
Thus, next to its main purpose of parallelization, a major contribution of \ug has been an improved stability of \scip.
For example, complete thread-safety of \scip was only achieved due to the development of \fscip.
 
Since \ug provides a systematic way to parallelize a state-of-the-art sequential or multi-threaded solver to run on a large scale distributed memory environment,
with version 1.0, \ug is generalized to a software framework for a {\em high-level task parallelization framework}\footnote{For concept of UG's high-level task parallelization framework, see \url{https://ug.zib.de/doc-1.0.0/html/CONCEPT.php}}.
That is, with version 1.0, \ug will not only parallelize the tree search of branch-and-bound based solvers, but allow the parallelization of other kind of solvers.
On top of that, \ug version 1.0 will also allow more flexibility and customization when parallelizing a branch-and-bound based solver for a specific purpose.
For an example, see the recent adaptation CMAP-LAP (Configurable Massively Parallel solver framework for LAttice Problems) of \ug to solve lattice problems~\cite{TateiwaShinanoYamamuraetal.2021}.

With the new beta version of \ug 1.0, which is released with the \scipoptv, \ug has caught up with interface changes in \scip and includes a few more bugfixes.
It does not include many new features. However, the possibility to appropriately specify an optimality gap limit has been added.

\subsection{Setting optimality gap limit}
For \fscip/\pscip, \scip parameters can be set by using command line options,``-sl, -sr, -s'' as below:
\begin{verbatim}
../bin/fscip 
syntax: ../bin/fscip fscip_param_file problem_file_name [-l <logfile>] [-q] 
      [-sl <settings>]  [-s <settings>] [-sr <root_settings>] [-w <prefix_warm>] 
      [-sth <number>] [-fsol <solution_file>] [-isol <initial solution file]
-l <logfile>           : copy output into log file
-q                     : suppress screen messages
-sl <settings>         : load parameter settings (.set) file for LC presolving
-s <settings>          : load parameter settings (.set) file for solvers
-sr <root_settings>    : load parameter settings (.set) file for root
-w <prefix_warm>       : warm start file prefix ( prefix_warm_nodes.gz and 
                         prefix_warm_solution.txt are read )
-sth <number>          : the number of solver threads used
-fsol <solution file> : specify output solution file
-isol <intial solution file> : specify initial solution file	
\end{verbatim}

Therefore, parallel solving algorithms of \fscip/\pscip can be controlled just by setting these parameters.
The optimality gap limit could be set by using the ``-s'' option. 
However, using default settings, \fscip/\pscip executes presolving in the LoadCoordinator,
which is the controller thread (process) of \fscip (\pscip, respectively),
and the presolved instance is passed on to all solvers~\cite{Shinano18fiber}. 
Therefore, the ``-s'' option applies to solving the presolved instance, not the original one. 
This can be a problem when trying to set a gap limit, for instance.
With a previous release, a \ug parameter was added to handle this appropriately, but it turned out to not work well and has now been removed again.
Instead, with this version, an optimality gap limit can be set for the original instance. 
To do so, the gap limit should be set in the \scip parameter setting file that is specified with the ``-sl'' option.
The LoadCoordinator will then handle the gap limit appropriately.

\section{The \gcg Decomposition Solver}
\label{sect:gcg}
\gcg is an extension that turns \scip into a branch-and-price or
branch-and-Benders-cut solver for mixed-integer linear programs. \gcg
can automatically detect a model structure that allows for a
Dantzig-Wolfe reformulation or Benders decomposition. The
reformulation process and the corresponding algorithmics like Benders
cut and column generation is done automatically without interaction from the users.
They just need to
provide the model. The latest version is \gcgv. Here are the few
changes since version 3.5 upon which we reported along with \scip
version 8.0.

\gcgv mainly contains code base improvements, with no major algorithmic
changes. Most importantly for developers, the API has mildly changed:
The prefix \texttt{DEC\_} was replaced with \texttt{GCG\_} to achieve
a consistent naming. 

The model structures detected by \gcg are called decompositions. 
The detection process itself was described in the SCIP 6.0 release
report. Part of this process is based on classifiers which group
constraints and variables for their potential roles in a decomposition.
From the usually many decompositions found, users can select manually or let
\gcg select based on different scores.  The score implementation was
refactored: Previously, it was cumbersome to add user-defined
scores to \gcg, but with \gcgv, new scores can be added as plugins.  Each
score must provide a function that calculates a score value for a
(partial) decomposition. The macro \texttt{GCG\_DECL\_SCORECALC} is
provided to declare the method that implements the scoring and is automatically called by \gcg.

The \texttt{display} dialog command has been extended and can now be used
to print information about registered scores and classifiers.
Furthermore, \gcgv now supports compiling with Microsoft Visual C++ (MSVC) and it is possible to use CMake to simplify the build process.

The Python interface \textsc{PyGCGOpt} mirrors the above changes.
Users can now customize the detection process even more by adding
own classifiers and scores.

\section{SCIP-SDP}
\label{sect:scip-sdp}

\scipsdp is a solver for handling mixed-integer semidefinite programs,
w.l.o.g, written in the following form
\begin{equation}\label{MISDP}
  \begin{aligned}
    \inf \quad & b^\T y \\
    \text{s.t.} \;\;\, & \sum_{k=1}^m A^k\, y_k - A^0 \succeq 0, \\
    & \ell_i \leq y_i \leq u_i && \forall\, i \in [m], \\
    & y_i \in \Z && \forall\, i \in I,
  \end{aligned}
\end{equation}
with symmetric matrices $A^k \in \R^{n \times n}$ for
$i \in \{0, \dots, m\}$, $b \in \R^m$, $\ell_i \in \R \cup \{- \infty\}$,
$u_i \in \R \cup \{\infty\}$ for all $i \in [m] \defi \{1, \dots,
m\}$. The set of indices of integer variables is given by $I \subseteq [m]$,
and $M \succeq 0$ denotes that a matrix $M$ is positive
semidefinite.

\scipsdp uses an SDP-based branch-and-bound approach based on \scip
(default). It also supports the possibility to use linear inequalities in
an LP-based approach. Which one is faster, depends on the instance.

The development of \scipsdp proceeded along a series of dissertations:
Mars~\cite{Mar13}, Gally~\cite{Gal19}, and
Matter~\cite{Mat22}. Corresponding articles
are~\cite{GallyPfetschUlbrich2018} (existence of Slater points in
branch-and-bound, dual fixing) and \cite{MatP22} (presolving
techniques). In the following, we give an overview of the main changes
since version 4.0.0, which was reported on in the \scip8 report. The
current version of \scipsdp is 4.3.0.

\subsection{Symmetry Handling of MISDPs}

One can handle permutation symmetries of~\eqref{MISDP} in the sense of
Section~\ref{sec:symmetry}. We will sketch how this works and refer to~\cite{HojP23} for details.

To define symmetries of~\eqref{MISDP}, for a matrix $A \in \R^{n \times n}$
and a permutation $\sigma$ of $[n]$, let
\[
  \sigma(A)_{ij} \defi A_{\sigma^{-1}(i),\sigma^{-1}(j)}\quad \forall\, i,\,j \in [n].
\]

\begin{definition}
  A permutation $\pi$ of variable indices $[m]$ is a \emph{formulation
    symmetry} of~\eqref{MISDP} if there exists a permutation $\sigma$ of
  the dimensions $[n]$ such that
  \begin{enumerate}
  \item $\pi(I) = I$, $\pi(\ell) = \ell$, $\pi(u) = u$, and $\pi(b) = b$\\
    ($\pi$ leaves integer variables, variable bounds, and the objective
    coefficients invariant),
  \item $\sigma(A^0) = A^0$ and, for all~$i \in [m]$, $\sigma(A^i) = A^{\pi^{-1}(i)}$.
  \end{enumerate}
\end{definition}

Such symmetries can be detected by using graph automorphism algorithms,
see~\cite{HojP23}. Examples of (formulation) symmetries computed for a
testset of \MISDP instances can also be found in~\cite{HojP23}.  Note that
we do not exploit symmetries in the matrix solutions of the SDPs like it
has been done in \cite{GatP04,HuSW22}, for example.

\begin{table}
  \caption{Results on a testset of 184 MISDP instances with/without using
    symmetry handling in \scip. We report the shifted geometric means of the
    running times in seconds and number of nodes. Column ``symtime'' and
    ``\# gens'' report the averge time for symmetry handling (including
    detection) and number of generators, respectively. The ``all optimal''
    block reports results for the 168 instances that were solved by both
    methods. The last column gives the shifted geometric mean running time
    only for the 21 instances that contain some symmetry.}
  \label{tab:MISDP-symmetries}
  \begin{tabular*}{\textwidth}{@{}l@{\extracolsep{\fill}}rrrrrr@{}}\toprule
    &\multicolumn{3}{c}{all (184) } & \multicolumn{2}{c}{all optimal (168)} & \multicolumn{1}{c@{}}{only symmetric (21)}\\\cmidrule{2-4}\cmidrule{5-6}\cmidrule{7-7}
    & time (s) & symtime (s) & \# gens & time (s) & \#nodes & time (s) \\\midrule
    without   &  130.6 & --   & -- & 95.0 & 778.3 & 45.07\\
    with      &  125.3 & 0.44 & 99 & 90.8 & 760.6 & 29.84\\\bottomrule
  \end{tabular*}
\end{table}

Table~\ref{tab:MISDP-symmetries} presents computational results -- we refer
to~\cite{HojP23} for the setup and details. We observe a speed-up of about
\SI{4}{\percent} for all instances and of about \SI{34}{\percent} for the
21 instances that contain symmetry.

\begin{figure}[t]
  \begingroup
  \definecolor{dgreen}{HTML}{228B22}
  \centering
  \begin{tikzpicture}
    [node/.style={circle,inner sep=2pt}]
    \node (d1) at (0,0) [node,fill=red,label=above:1] {};
    \node (d2) at (3,0) [node,fill=red,label=above:2] {};
    \node (d3) at (6,0) [node,fill=red,label=above:3] {};

    \node (e1) at (0.5,-1.5) [node,fill=orange,label={[label distance=-1ex]left:{\tiny$(1,1)^1$}}] {};
    \node (e2) at (1.5,-1.5) [node,fill=blue,label={[label distance=-1ex]left:{\tiny$(1,2)^1$}}] {};
    \node (e3) at (2.5,-1.5) [node,fill=blue,label={[label distance=-1ex]left:{\tiny$(2,2)^1$}}] {};

    \node (f1) at (3.5,-1.5) [node,fill=blue,label={[label distance=-1ex]right:{\tiny$(2,2)^2$}}] {};
    \node (f2) at (4.5,-1.5) [node,fill=blue,label={[label distance=-1ex]right:{\tiny$(2,3)^2$}}] {};
    \node (f3) at (5.5,-1.5) [node,fill=orange,label={[label distance=-1ex]right:{\tiny$(3,3)^2$}}] {};

    \node (c1) at (7.0,-1.5) [node,fill=blue,label={[label distance=-1ex]below:{\tiny$(1,1)^0$}}] {};
    \node (c2) at (8.0,-1.5) [node,fill=blue,label={[label distance=-1ex]below:{\tiny$(2,2)^0$}}] {};
    \node (c3) at (9.0,-1.5) [node,fill=blue,label={[label distance=-1ex]below:{\tiny$(3,3)^0$}}] {};
    
    \node (v1) at (1.5,-3) [node,fill=dgreen,label=left:{$y_1$}] {};
    \node (v2) at (4.5,-3) [node,fill=dgreen,label=left:{$y_2$}] {};

    \draw[-] (d1) -- (e1);
    \draw[-] (d1) -- (e2);
    \draw[-] (d2) -- (e2);
    \draw[-] (d2) -- (e3);

    \draw[-] (d2) -- (f1);
    \draw[-] (d2) -- (f2);
    \draw[-] (d3) -- (f2);
    \draw[-] (d3) -- (f3);

    \draw[-] (d1) -- (c1);
    \draw[-] (d2) -- (c2);
    \draw[-] (d3) -- (c3);

    \draw[-] (e1) -- (v1);
    \draw[-] (e2) -- (v1);
    \draw[-] (e3) -- (v1);

    \draw[-] (f1) -- (v2);
    \draw[-] (f2) -- (v2);
    \draw[-] (f3) -- (v2);
  \end{tikzpicture}
  \caption{Illustration of symmetry detection graph.}
  \label{fig:detection}
  \endgroup
\end{figure}
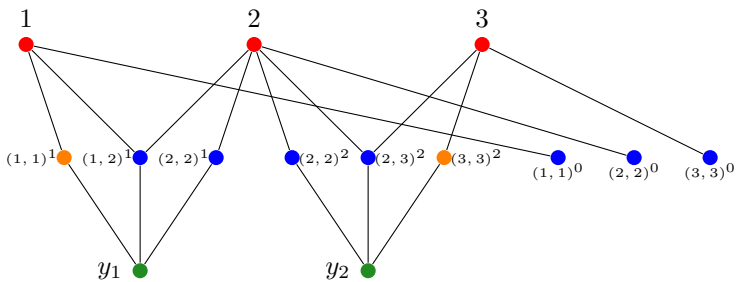

Since the appearance of~\cite{HojP23}, \scipsdp has been changed to allow
for using the callback access to symmetry computation in \scip, see
Section~\ref{sec:symdetect}. Thereby, we have also changed the symmetry
detection graph. It now only contains a single node for each nonzero entry
of the matrix that is connected to the `dimension'-nodes. The corresponding
graph for the following \MISDP is given in Figure~\ref{fig:detectsym}.
\[
  \inf\;\left\{
    y_1 + y_2
    \st
    \left(
      \begin{smallmatrix}
        3 & 1 & 0\\
        1 & 1 & 0\\
        0 & 0 & 0
      \end{smallmatrix}
    \right)
    y_1
    +
    \left(
      \begin{smallmatrix}
        0 & 0 & 0\\
        0 & 1 & 1\\
        0 & 1 & 3
      \end{smallmatrix}
    \right)
    y_2
    -
    \left(
      \begin{smallmatrix}
        1 & 0 & 0\\
        0 & 1 & 0\\
        0 & 0 & 1
      \end{smallmatrix}
    \right)
    \succeq 0,\;
    0 \leq y_1,\, y_2 \leq 1,\;
    y_1,\, y_2 \in \Z
  \right\}.
\]

In Figure~\ref{fig:detection}, the colors of the nodes allow to distinguish
different types and values. The topmost nodes represent the
dimensions. Node $(i,j)^k$ represents the symmetric entries $(i,j)$ and
$(j,i)$ (or diagonal entry $(i,i)$) of matrix $A^k$. The colors of these
nodes correspond to distinct coefficients in the matrices. The only
non-trivial color-preserving automorphism of the graph
exchanges~$y_1 \leftrightarrow y_2$, ${(1,1)^1 \leftrightarrow (3,3)^2}$,
${(1,2)^1 \leftrightarrow (2,3)^2}$, ${(2,2)^1 \leftrightarrow (2,2)^2}$,
${(1,1)^0 \leftrightarrow (3,3)^0}$,
$1 \leftrightarrow 3$, and keeps node~2 fixed.  This leads to the variable
permutation~$\pi$, which exchanges~$y_1$ and~$y_2$, and the matrix
permutation~$\sigma$, which exchanges~$1$ and~$3$.

\subsection{Conflict Analysis for \MISDPs}

The original idea of conflict analysis was to learn from infeasible nodes
in a branch-and-bound-tree. To this end, Achterberg~\cite{Ach07} transferred
ideas from SAT-solving to \MILPs.
One further way is to try to learn
cuts from solutions of the duals, which is called ``dual ray/solution
analysis'' in Witzig et al.~\cite{WitBH17} and Witzig~\cite{Wit22}.

To briefly explain the application to \MISDPs, consider the SDP relaxation
of~\eqref{MISDP}. Given a positive semidefinite $\hat{X} \in \R^{n \times n}$,
we observe that the inner product with a positive semidefinite matrix $M
\in \R^{n \times n}$
is nonnegative:
\[
  \iprod{\hat{X}}{M} \defi \sum_{i,j = 1}^n \hat{X}_{ij}\, M_{ij} \geq 0.
\]
Thus, defining $A(y) \defi \sum_{k=1}^m A^k\, y_i$, we get
\begin{equation}\label{eq:MISDPConflictInequality}
  \iprod{\hat{X}}{A(y)} = \sum_{k=1}^m \iprod{\hat{X}}{A^k}\, y_k \geq 0
\end{equation}
for every feasible solution~$y$ of~\eqref{MISDP}. Note that this is a
(redundant) linear inequality in~$y$. The idea is to use it in the
propagation of variable bounds and not explicitly add it to~\eqref{MISDP}.

There are two natural ways to obtain good candidates for $\hat{X}$. If the
relaxation is feasible, we obtain a solution
$(\hat{X}, \hat{r}^\ell, \hat{r}^u)$ of the dual
\[
  \begin{aligned}
    \sup \quad & \iprod{A^0}{X} + \ell^\T r^\ell - u^\T r^u\\
    \text{s.t.} \quad & \iprod{A^j}{X} + r^\ell_j - r^u_j = b_j && \forall \, j \in [m],\\
    & X \succeq 0,\; r^\ell,\; r^u \geq 0.
  \end{aligned}
\]
Similarly, if the relaxation is infeasible and a constraint qualification
holds, one can obtain a dual ray satisfying
\[
  \begin{aligned}
    & \iprod{A^j}{X} + r^\ell_j - r^u_j = 0 && \forall\; j \in [m],\\
    & \iprod{A^0}{X} + \ell^\T r^\ell - u^\T r^u > 0,\\
    & X \succeq 0,\; r^{\ell},\; r^u \geq 0.
  \end{aligned}
\]
One can prove that~\eqref{eq:MISDPConflictInequality} is infeasible
with respect to the local bounds $\ell$ and $u$ and can therefore provide
a proof of infeasibility, see~\cite{Pfe23}.

\scipsdp generates a conflict constraint~\eqref{eq:MISDPConflictInequality}
for each feasible or infeasible node, stores them as constraints, and
performs bound propagation. This leads to a speed-up and node reduction of
about 8\% on the same testset used in the previous section. We refer
to~\cite{Pfe23} for more details.

\subsection{Further Changes}

Similar to \scip, the license of \scipsdp has changed to Apache 2.0.

Several improvements have been made to speed up some of the presolving
methods. One can use ARPACK instead of Lapack for eigenvalue computations
(use \texttt{ARPACK = true} when using makefiles); this is usually slower
for the typical sizes of \MISDPs. Moreover, when running the \LP-based
approach (\texttt{misc/solvesdps = 0}), now CMIR inequalities are generated
by default. To implement conflict analysis, the handling of dual
solutions has been extended and improved.

\section{Final Remarks}
\label{sect:finalRemarks}

The \scipoptv release provides new functionality along with improved performance and
reliability.
In \scip, the changes to the symmetry detection feature include new techniques for handling symmetries of non-binary variables, restructuring of the mechanism to detect symmetries of the custom constraints, and detection of the signed permutation symmetries.
New interfaces to \nauty~\cite{Nauty} as well as the preprocessor \sassy~\cite{Sassy} have been added.
A new nonlinear handler for signomial functions and improvements to the existing nonlinear handler for quadratic expressions were implemented.
A new diving heuristic for handling indicator constraints that are used to represent the semi-continuous variables, an extension of the dynamic partition search heuristic,
and a new adaptive heuristic that dynamically adapts the application of large neighborhood search and diving heuristics to the characteristics of the current instance were also implemented.
Furthermore, a new separator called the Lagromory separator for generating potentially lower-ranked cuts and reducing the dual bound stalling due to the dual degeneracy, and two new cut selection schemes were included: \texttt{cutsel/ensemble} that adapts with respect to the given instance properties and \texttt{cutsel/dynamic} that aims to enhance the near-orthogonal threshold methodology used in the default cut selection scheme.
A new branching criterion called the GMI branching was implemented. It is available both as a stand-alone rule and also integrated within the default hybrid branching rule of \scip. It considers a new scoring component based on the GMI cuts corresponding to the fractional variable in a given \LP solution.
Finally, a new interface to the \highs~\LP solver along with technical improvements to the AMPL reader and OBBT propagator were implemented.

Regarding usability, various interfaces were improved and new interfaces were added.
The AMPL interface to \scip was extended to support parameters from the AMPL command scripts.
The Java interface to \scip, \jscipopt, is being maintained actively and was extended with new functionality.
The Python interface to \scip, \pyscipopt, can now be fully installed using PyPI.
A new Python package, \pyscipoptml, is available to automatically formulate machine learning models into MIPs.
The Julia interface to \scip, \scipjl, was extended to be able to access additional plugins of \scip.
Two new interfaces to \scip, Rust interface called \russcip and C++ interface called \scippp, are also available, along with a new Python interface to \soplex called \pysoplex.

The \LP solver \soplex now supports incremental precision boosting for exact \LP solving over the rational numbers. It is available as a stand-alone as well as in combination with the existing \LP iterative refinement approach.
The presolving library \papilo now has a new feature called proof logging that allows the generation of machine-verifiable certificates for presolving of binary problems to be able to prove the correctness of the computations.
The parallel framework \ug now has a new beta version of \ug~1.0 that includes the latest interface changes of \scip along with new bugfixes.
A new feature to appropriately set the gap limit has also been added.
The \gcg decomposition solver now includes improvements to its code base, can be compiled with Microsoft Visual C++ (MSVC), and supports CMake as a build system.
The \scip extension \scipsdp has been improved to include new symmetry handling techniques and conflict analysis for \MISDPs, along with other improvements in its presolving methods and cut generation techniques.

These developments yield an overall performance improvement of both \MILP and \MINLP benchmarking instances.
\scipv is able to solve~19 more \MILP instances as compared to \scip~8.0 with a speedup of~2\% on the affected instances.
This speedup further increases to~6\% when only the hard instances requiring at least 1000 seconds by at least one setting are considered.
The number of nodes required for \MILPs that were solved by both versions of \scip also reduce considerably by~17\% in \scipv.
These performance gains are more prominent for \MINLPs.
\scipv solves~5 more \MINLP instances as compared to \scip~8.0 with performance improvements of~4\% in time (for all \MINLPs) and~13\% in the number of nodes (for the \MINLPs that were solved by both the versions of \scip).
When looking at hard instances requiring at least 1000 seconds by at least one setting, the gains are further increased to~20\% and~46\% in time and number of nodes, respectively.
Furthermore, when restrcited to nonconvex instances only, \scipv is faster by~8\%.
Hence, \scipv has become faster and more reliable as compared to \scip~8.0.
\subsection*{Acknowledgements}

The authors want to thank all previous developers and contributors to the \scipopt and
all users that reported bugs and often also helped reproducing and fixing the bugs.
Thanks to Herman Appelgren, Gerald Gamrath, Stephen J. Maher, Daniel Rehfeldt, and
Michael Winkler, for various contributions.
Thanks to Julian Hall for supporting the creation of the \highs~\LP interface.

\subsection*{Contributions of the Authors}
The material presented in the article is highly related to code and software.
In the following, we try to make the corresponding contributions of the authors and possible contact points more transparent.

\begin{itemize}
 \item JvD implemented the generalization of symmetry handling methods for
   binary variables to general variables (\Cref{sec:symmetryhandling}).
 \item CH implemented the symmetry detection framework via symmetry detection
   graph callbacks in constraint handlers (\Cref{sec:symdetect}).
 \item MP implemented the interface to \nauty, and the interface to allow \sassy as preprocessor for \bliss (\Cref{sec:symmetryinterfaces}).
 \item MP and GM implemented the interface to allow \sassy as preprocessor for \nauty (\Cref{sec:symmetryinterfaces}).
 \item LX implemented the nonlinear handler for signomial functions (\Cref{subsubsect:signomial}).
 \item AC and FeS implemented monoidal strengthening for quadratic constraints (\Cref{subsubsect:quadratic}).
 \item KH and AH implemented the indicator diving heuristic (\Cref{subsubsect:indicatordiving}).
 \item KH extended the DPS heuristic (\Cref{subsubsect:dps}).
 \item AC implemented the online scheduling procedure of primal heuristics (\Cref{subsubsect:adaptiveheuristics}).
 \item SB implemented the Lagromory separator (\Cref{subsubsect:cutgeneration}).
 \item MT and MB worked on the ensemble cut selector (\Cref{subsubsect:cutselection}).
 \item CG worked on cut statistics and dynamic cut selection (\Cref{subsubsect:cutselection}).
 \item MT and MB worked on the GMI branching and hybrid branching improvements (\Cref{subsubsec:GMI} and~\Cref{subsubsec:GMI_improve}).
 \item AG, GM, and AH developed the \highs~\LP interface (\Cref{subsect:highs}).
 \item SV extended the AMPL \texttt{.nl} reader (\Cref{subsect:further}).
 \item KH extended the OBBT propagator (\Cref{subsect:further}).
 \item LE worked on the \soplex solver (\Cref{sect:soplex}).
 \item AH implemented the proof-logging in \papilo (\Cref{sect:papilo}).
 \item SV extended the AMPL interface (\Cref{subsect:ampl}).
 \item KK is the current maintainer of \jscipopt and has worked on the developments listed in
\Cref{subsect:jscipopt} (other than the ones where different developer names are mentioned).
 \item JD actively maintained and contributed to \pyscipopt (\Cref{subsect:pyscipopt}).
 \item MT worked on the pip installation of \pyscipopt through PyPI (\Cref{subsect:pyscipopt}).
 \item MB worked on the Julia interface \scipjl (\Cref{subsect:scipjl}).
 \item MG implemented the \russcip interface (\Cref{subsect:russcip}).
 \item IH programmed the initial release of \scippp (\Cref{subsect:scippp}).
 \item SB implemented the \pysoplex interface (\Cref{subsect:pysoplex}).
 \item YS worked on the \ug framework (\Cref{sect:ug})
 \item Concerning \gcg (\Cref{sect:gcg}), EM and JL refactored and improved the code (general API changes and scores) and added the MSVC support; JL and TD extended \solver{PyGCGOpt}.
 \item CH and MP implemented symmetry handling for \scipsdp and MP implemented
   conflict analysis in \scipsdp (\Cref{sect:scip-sdp}).
 \item RvdH contributed to various bug fixes and file reader updates.
 \item DK fixed several bugs and contributed to other fixes regarding reliability.
   Furthermore, DK implemented an algorithm to simplify the debugging process by generating a potentially smaller-sized instance and settings that can reproduce the bugs faster.
 \item JM and FrS worked on the continuous integration system, regular testing, binary distributions, and website development. JM is the contact person for these aspects.
\end{itemize}

\renewcommand{\refname}{\normalsize References}
\setlength{\bibsep}{0.25ex plus 0.3ex}
\bibliographystyle{abbrvnat}

\begin{small}
\bibliography{scipopt}

\end{small}

\subsection*{Author Affiliations}
\hypersetup{urlcolor=black}
\newcommand{\myorcid}[1]{ORCID: \href{https://orcid.org/#1}{#1}}
\newcommand{\myemail}[1]{E-mail: \href{#1}{#1}}
\newcommand{\myaffil}[2]{{\noindent #1}\\{#2}\bigskip}

\small

\myaffil{Suresh Bolusani}{
  Zuse Institute Berlin, Department AIS$^2$T, Takustr.~7, 14195~Berlin, Germany\\
  \myemail{bolusani@zib.de}\\
  \myorcid{0000-0002-5735-3443}}

\myaffil{Mathieu Besançon}{%
  Université Grenoble Alpes, Inria, LIG, 38000 Grenoble, France, and\\
  Zuse Institute Berlin, Department AIS$^2$T, Takustr.~7, 14195~Berlin, Germany\\
  \myemail{besancon@zib.de}\\
  \protect\myorcid{0000-0002-6284-3033}}

\myaffil{Ksenia Bestuzheva}{
  Zuse Institute Berlin, Department AIS$^2$T, Takustr.~7, 14195~Berlin, Germany\\
  \myemail{bestuzheva@zib.de}\\
  \myorcid{0000-0002-7018-7099}}

\myaffil{Antonia Chmiela}{%
  Zuse Institute Berlin, Department AIS$^2$T, Takustr.~7, 14195~Berlin, Germany\\
  \myemail{chmiela@zib.de}\\
  \myorcid{0000-0002-4809-2958}}

\myaffil{Jo{\~{a}}o Dion{\'{i}}sio}{%
CMUP and Department of Computer Science, Faculty of Sciences, University of Porto, R Campo Alegre, 4169–007 Porto, Portugal\\
  \myemail{joao.goncalves.dionisio@gmail.com}\\
  \myorcid{0009-0005-5160-0203}}

\myaffil{Tim Donkiewicz}{%
  RWTH Aachen University, Lehrstuhl f\"ur Operations Research, Kackertstr.~7, 52072~Aachen, Germany\\
  \myemail{tim.donkiewicz@rwth-aachen.de}\\
  \myorcid{0000-0002-5721-3563}}

\myaffil{Jasper van Doornmalen}{%
  Eindhoven University of Technology, Department of Mathematics and Computer Science, P.O.\ Box 513, 5600 MB Eindhoven, The Netherlands\\
  \myemail{m.j.v.doornmalen@tue.nl}\\
  \myorcid{0000-0002-2494-0705}}

\myaffil{Leon Eifler}{%
  Zuse Institute Berlin, Department AIS$^2$T, Takustr.~7, 14195~Berlin, Germany\\
  \myemail{eifler@zib.de}\\
  \myorcid{0000-0003-0245-9344}}

\myaffil{Oliver Gaul}{%
  RWTH Aachen University, Lehrstuhl f\"ur Operations Research, Kackertstr.~7, 52072~Aachen, Germany\\
  \myemail{oliver.gaul@rwth-aachen.de}\\
  \myorcid{0000-0002-2131-1911}}

\myaffil{Mohammed Ghannam}{%
  Zuse Institute Berlin, Department AIS$^2$T, Takustr.~7, 14195~Berlin, Germany\\
  \myemail{ghannam@zib.de}\\
  \myorcid{0000-0001-9422-7916}}

\myaffil{Ambros Gleixner}{%
  Zuse Institute Berlin, Department AIS$^2$T, Takustr.~7, 14195~Berlin, Germany\\
  \myemail{gleixner@zib.de}\\
  \myorcid{0000-0003-0391-5903}}

\myaffil{Christoph Graczyk}{%
  Zuse Institute Berlin, Department AIS$^2$T, Takustr.~7, 14195~Berlin, Germany\\
  \myemail{graczyk@zib.de}\\
  \myorcid{0000-0001-8990-9912}}

\myaffil{Katrin Halbig}{%
  Friedrich-Alexander Universität Erlangen-Nürnberg, Department of Data Science, Cauerstr.~11, 91058~Erlangen, Germany\\
  \myemail{katrin.halbig@fau.de}\\
  \myorcid{0000-0002-8730-3447}}

\myaffil{Ivo Hedtke}{%
  Schenker AG, Global Data \& AI, Kruppstr.~4, 45128~Essen, Germany\\
  \myemail{ivo.hedtke@dbschenker.com}\\
  \myorcid{0000-0003-0335-7825}}

\myaffil{Alexander Hoen}{%
  Zuse Institute Berlin, Department AIS$^2$T, Takustr.~7, 14195~Berlin, Germany\\
  \myemail{hoen@zib.de}\\
  \myorcid{0000-0003-1065-1651}}

\myaffil{Christopher Hojny}{%
  Eindhoven University of Technology, Department of Mathematics and Computer Science, P.O.\ Box 513, 5600 MB Eindhoven, The Netherlands\\
  \myemail{c.hojny@tue.nl}\\
  \myorcid{0000-0002-5324-8996}}

\myaffil{Rolf van der Hulst}{%
  University of Twente, Department of Discrete Mathematics and Mathematical Programming, P.O. Box 217, 7500 AE Enschede, The Netherlands\\
  \myemail{r.p.vanderhulst@utwente.nl}\\
  \myorcid{0000-0002-5941-3016}}

\myaffil{Dominik Kamp}{%
  University of Bayreuth, Chair of Economathematics, Universitaetsstr.~30, 95440~Bayreuth, Germany\\
  \myemail{dominik.kamp@uni-bayreuth.de}\\
  \myorcid{0009-0005-5577-9992}}

\myaffil{Thorsten Koch}{%
  Technische Universit\"at Berlin, Chair of Software and Algorithms for Discrete Optimization, Stra\ss{}e des 17. Juni 135, 10623 Berlin, Germany, and\\
  Zuse Institute Berlin, Department A$^2$IM, Takustr. 7, 14195~Berlin, Germany\\
  \myemail{koch@zib.de}\\
  \myorcid{0000-0002-1967-0077}}

\myaffil{Kevin Kofler}{%
  DAGOPT Optimization Technologies GmbH\\
  \myemail{kofler@dagopt.com}}

\myaffil{Jurgen Lentz}{%
  RWTH Aachen University, Lehrstuhl f\"ur Operations Research, Kackertstr.~7, 52072~Aachen, Germany\\
  \myemail{jurgen.lentz@rwth-aachen.de}\\
  \myorcid{0009-0000-0531-412X}}

\myaffil{Julian Manns}{%
  Zuse Institute Berlin, Department AIS$^2$T, Takustr.~7, 14195~Berlin, Germany\\
  \myemail{manns@zib.de}}

\myaffil{Gioni Mexi}{%
  Zuse Institute Berlin, Department AIS$^2$T, Takustr.~7, 14195~Berlin, Germany\\
  \myemail{mexi@zib.de}\\
  \myorcid{0000-0003-0964-9802}}

\myaffil{Erik M\"uhmer}{%
  RWTH Aachen University, Lehrstuhl f\"ur Operations Research, Kackertstr.~7, 52072 Aachen, Germany\\
  \myemail{erik.muehmer@rwth-aachen.de}\\
  \myorcid{0000-0003-1114-3800}} 

\myaffil{Marc E.~Pfetsch}{%
  Technische Universität Darmstadt, Fachbereich Mathematik, Dolivostr.~15, 64293~Darmstadt, Germany\\
  \myemail{pfetsch@mathematik.tu-darmstadt.de}\\
  \myorcid{0000-0002-0947-7193}}

\myaffil{Franziska Schlösser}{%
  Fair Isaac Germany GmbH, Takustr.~7, 14195~Berlin, Germany\\
  \myemail{franziskaschloesser@fico.com}}

\myaffil{Felipe Serrano}{%
  COPT GmbH, Berlin, Germany\\
  \myemail{serrano@copt.de}\\
  \myorcid{0000-0002-7892-3951}}

\myaffil{Yuji Shinano}{%
  Zuse Institute Berlin, Department A$^2$IM, Takustr.~7, 14195~Berlin, Germany\\
  \myemail{shinano@zib.de}\\
  \myorcid{0000-0002-2902-882X}}

\myaffil{Mark Turner}{%
  Zuse Institute Berlin, Department A$^2$IM, Takustr.~7, 14195~Berlin, Germany\\
  \myemail{turner@zib.de}\\
  \myorcid{0000-0001-7270-1496}}

\myaffil{Stefan Vigerske}{%
  GAMS Software GmbH, c/o Zuse Institute Berlin, Department AIS$^2$T, Takustr.~7, 14195~Berlin, Germany\\
  \myemail{svigerske@gams.com}\\
  \myorcid{0009-0001-2262-0601}}

\myaffil{Dieter Weninger}{%
  Friedrich-Alexander Universität Erlangen-Nürnberg, Department of Mathematics, Cauerstr.~11, 91058~Erlangen, Germany\\
  \myemail{dieter.weninger@fau.de}\\
  \myorcid{0000-0002-1333-8591}}
  
\myaffil{Liding Xu}{%
  École polytechnique, LIX CNRR,  Rue Honoré d'Estienne d'Orves.~1, 9120~Palaiseau, France\\
  \myemail{lidingxu.ac@gmail.com}\\
  \myorcid{0000-0002-0286-1109}}

\end{document}